\documentclass[oneside,english]{siamltex}
\usepackage[T1]{fontenc}
\usepackage[latin9]{inputenc}
\usepackage{geometry}
\usepackage{xcolor}
\usepackage{pdfcolmk}
\usepackage{array}
\usepackage{multirow}
\usepackage{amsbsy}
\usepackage{amssymb}
\usepackage{graphicx}
\usepackage[numbers]{natbib}
\PassOptionsToPackage{normalem}{ulem}
\usepackage{ulem}

\makeatletter

\providecommand{\tabularnewline}{\\}
\providecolor{lyxadded}{rgb}{0,0,1}
\providecolor{lyxdeleted}{rgb}{1,0,0}

\DeclareRobustCommand{\lyxsout}[1]{\ifx\\#1\else\sout{#1}\fi}

\renewcommand\eqref[1]{(\ref{#1})}

\usepackage[english]{babel}
\renewcommand\eqref[1]{(\ref{#1})}

\@ifundefined{showcaptionsetup}{}{%
 \PassOptionsToPackage{caption=false}{subfig}}
\usepackage{subfig}
\makeatother

\usepackage{babel}
\begin{document}

\title{A Posteriori Error Analysis of Fluid-Stucture Interactions: Time
Dependent Error}

\author{Jay A. Stotsky\thanks{University of Colorado, Department of Applied Mathematics (\email{jay.stotsky@colorado.edu, dmbortz@colorado.edu})}
\and 
David M. Bortz\footnotemark[1]}
\maketitle
\begin{abstract}
\emph{A posteriori} error analysis is a technique to quantify the error in particular simulations of a numerical approximation method. 
In this article, we use such an approach to analyze
how various error components propagate in certain moving boundary
problems. We study quasi-steady state simulations where slowly moving
boundaries remain in mechanical equilibrium with a surrounding fluid.
Such problems can be numerically approximated with the \emph{Method
of Regularized Stokelets }(MRS), a popular method used for studying
viscous fluid-structure interactions, especially in biological applications.
Our approach to monitoring the regularization error of the MRS is
novel, along with the derivation of linearized adjoint equations to
the governing equations of the MRS with a elastic elements. Our main
numerical results provide a clear illustration of how the error evolves
over time in several MRS simulations.
\end{abstract}

\section{Introduction}

The method of regularized Stokeslets (MRS) is a commonly used method
for studying viscous flow phenomena. The method is grid-free and based
on the numerical discretization of integrodifferential equations involving
singularity solutions of the Stokes equation. In this paper we develop
and apply \emph{a posteriori} error estimation techniques described
in \citep{collins_posteriori_2015} to the ordinary differential equations
that arise from the spatial discretization of the singular integral
operators associated with the MRS. Such systems frequently arise in
the study of quasi-steady state phenomena in biofluids \citep{aranda2015model,cortez_accuracy_1998,cortez_method_2001,cortez_method_2005,wrobel_enhanced_2016,wrobel_modeling_2014}.
However, little work has been applied towards understanding the error
in specific simulations that use the MRS (although see \citep{cortez_method_2005}
for some examples of numerical convergence). The main contributions
of this paper are the derivation of adjoints for the quasi-steady
state MRS and computational results based on a decomposition of the
numerical error into terms specifically tied each step in the discretization.
Furthermore, our approach towards computing the regularization error
during MRS simulations by treating the regularized equations as a
perturbation to a set of non-regularized equations is novel.

As explained in \citep{collins_posteriori_2015}, we see through
\emph{a posteriori }analysis that numerical error is composed of quadrature
error, extrapolation (or explicit) error, residual error, and regularization
error. The residual and regularization errors are due to the use of
a finite dimensional function space where the solution is approximated
and the use of regularization to remove singularities from the governing
equations. On the other hand, the quadrature and extrapolation errors
are due to the choice of numerical method. Quadrature arises from
the numerical approximation of integrals over each time interval.
For explicit methods, extrapolations of the numerical solution across
each time interval introduce error.  Since the standard time-stepping
methods for the MRS are explicit Runge-Kutta methods, we use the \emph{nodally
equivalent finite element method} formalism introduced in \citep{collins_posteriori_2015}
to derive polynomial solutions correspond to the Runge-Kutta solutions
over each time step in the solution. 

Section \ref{sec:Review-of-the} provides a brief overview of the
method of regularized Stokeslets. For readers familiar with the MRS,
it may be skimmed with the exception of Table \ref{tab:Notation}
which contains essential notation.Similarly, readers familiar with
\emph{a posteriori }error estimation techniques applied to systems
of ODEs may skim Section \ref{sec:A-Posteriori-Error}.

\section{Review of the Method of Regularized Stokeslets \label{sec:Review-of-the}}

In this section we briefly overview the the singular integral formulation
of Stokes flow.  The notation used in the rest of this section and
the proceeding sections is listed in Table \ref{tab:Notation}.

\begin{table}
\begin{centering}
\begin{tabular}{|c|l|}
\hline 
Symbol & Definition\tabularnewline
\hline 
\hline 
$t$ & time\tabularnewline
\hline 
$\boldsymbol{u}$ & velocity\tabularnewline
\hline 
$P$ & pressure\tabularnewline
\hline 
$\boldsymbol{\alpha}$ & Lagrangian coordinate\tabularnewline
\hline 
$k$ & discrete Lagrangian index\tabularnewline
\hline 
$r$ & radial distance, $r=\sqrt{x_{i}x_{i}}$.\tabularnewline
\hline 
$\zeta_{\delta}(r)$ & radially symmetric regularization kernel\tabularnewline
\hline 
$\boldsymbol{U}(\boldsymbol{x})$ & fundamental solution of Stokes equation\tabularnewline
\hline 
$\boldsymbol{U}_{\delta}(\boldsymbol{x})$ & smooth regularization of $\boldsymbol{U}(\boldsymbol{x})$ defined
by $\boldsymbol{U}_{\delta}(\boldsymbol{x})=\int\zeta_{\delta}(\boldsymbol{y})\boldsymbol{U}(\boldsymbol{x}-\boldsymbol{x})d\boldsymbol{y}$\tabularnewline
\hline 
$\boldsymbol{x}(\boldsymbol{\alpha},t)$ & spatially continuous flow map\tabularnewline
\hline 
$\boldsymbol{f}(\boldsymbol{x},t)$ & force\tabularnewline
\hline 
$\boldsymbol{x}_{k}(t)$ & spatially discrete, time-continuous flow map \tabularnewline
\hline 
$\boldsymbol{X}_{k}(t)$ & spatially discrete numerical approximation of $\boldsymbol{x}_{k}(t)$\tabularnewline
\hline 
$\boldsymbol{e}_{k}(t)$ & numerical error: $\boldsymbol{x}_{k}(t)-\boldsymbol{X}_{k}(t)$\tabularnewline
\hline 
$\left(\cdot,\cdot\right)$ & finite dimensional inner product e.g. $\left(f,g\right)=\sum_{k}f_{k}g_{k}$\tabularnewline
\hline 
$\left|\cdot\right|$ & norm over a finite dimensional vector space $\left|x\right|=(x,x)^{1/2}$\tabularnewline
\hline 
$\left\langle \cdot,\cdot\right\rangle $ & Inner product over time defined as $\int_{0}^{t}\left(\boldsymbol{u},\boldsymbol{v}\right)dt$\tabularnewline
\hline 
$\Vert\cdot\Vert$ & Norm over time and space\tabularnewline
\hline 
\multirow{2}{*}{$\left\langle \cdot,\cdot\right\rangle _{d}$} & discrete inner product approximation through a quadrature rule\tabularnewline
 & $\left\langle f,g\right\rangle _{d}=\sum_{q}^{Q}w_{q}f(t_{q})g(t_{q})$
with $\left\{ w_{q},t_{q}\right\} $ specified by the quadrature rule\tabularnewline
\hline 
$\mathcal{V}^{q}(I)$ & space of $q$th degree polynomials over an interval $I=[a,b]$. \tabularnewline
\hline 
$\mathcal{V}^{q}([0,T])$ & space of piecewise $q$th degree polynomials on $[0,T]$. \tabularnewline
\hline 
$\mathcal{P}$ & a projection operator related to Runge Kutta methods \tabularnewline
\hline 
$\pi$ & projection operator from a function space $\mathcal{V}$ to $\mathcal{V}^{q}(I)$\tabularnewline
\hline 
$A^{\ast}$ & adjoint, defined by $\left\langle A\boldsymbol{x},\boldsymbol{y}\right\rangle =\left\langle \boldsymbol{x},A^{\ast}\boldsymbol{y}\right\rangle $\tabularnewline
\hline 
$\bar{A}$ & trajectory-averaged operator, $\bar{A}=\int_{0}^{1}A(s\boldsymbol{x}+(1-s)\boldsymbol{y})ds$\tabularnewline
\hline 
\end{tabular}
\par\end{centering}
\caption{\label{tab:Notation}Definitions of frequently used terms.. }
\end{table}

\subsection{Stokes Equations and the Method of Regularized Stokeslets}

Stokes equations are an approximation of the Navier-Stokes equations
valid for viscous fluids at small length scales and low velocities.
Being an approximation method for solving Stokes equations, the Method
of Regularized Stokeslets, is often applied to problems in biofluids
which typically satisfy Stokes flow conditions. 

 The partial differential equations governing Stokes flow can be
written as
\begin{eqnarray*}
\nabla P & = & \Delta\boldsymbol{u}+\boldsymbol{f}\\
\nabla\cdot\boldsymbol{u} & = & 0.
\end{eqnarray*}
There exist fundamental velocity and pressure solutions, $U_{ij}(\boldsymbol{x}-\boldsymbol{y})$
and $P_{i}(\boldsymbol{x}-\boldsymbol{y})$, that solve the following
equations written in index notation\footnote{For a vector quantity, $u_{i,k}\equiv\frac{\partial u_{i}}{\partial x_{k}}$.
Summation over repeated indices is implied, i.e. in $\mathbb{R}^{3}$,
$u_{i,kk}\equiv\frac{\partial^{2}u_{i}}{\partial x_{1}^{2}}+\frac{\partial^{2}u_{i}}{\partial x_{2}^{2}}+\frac{\partial^{2}u_{i}}{\partial x_{3}^{2}}$}
\begin{eqnarray*}
P_{k,i}(\boldsymbol{x}-\boldsymbol{y}) & = & U_{ki,jj}(\boldsymbol{x}-\boldsymbol{y})+\delta(\boldsymbol{x}-\boldsymbol{y})\hat{e}_{i}^{(k)}\\
U_{kj,j}(\boldsymbol{x}-\boldsymbol{y}) & = & 0
\end{eqnarray*}
The basis vector, $\hat{\boldsymbol{e}}^{(k)}$ with $k\in\left\{ 1,2,3\right\} $
is defined such that $\hat{\boldsymbol{e}}^{(1)}$ is understood to
be the basis vector that points along the $x$-axis, $\hat{\boldsymbol{e}}^{(2)}$
along the $y$-axis, and $\hat{\boldsymbol{e}}^{(3)}$ along the $z$-axis.
The symbol, $\delta(\cdot)$ represents the Dirac delta distribution.
 Due to the linearity of Stokes equations, the solution for an arbitrary
vector, $\boldsymbol{g}\in\mathbb{R}^{3}$, is $u_{i}(\boldsymbol{r})=U_{ij}(\boldsymbol{r})g_{j}$.
The tensor-valued fundamental velocity solution denoted, $U_{ij}$
or $\boldsymbol{U}$ depending on the context, can be written in two
dimensions as  
\[
U_{ij}(\boldsymbol{r})=\frac{1}{4\pi}\left(-\delta_{ij}\log r+\frac{r_{i}r_{j}}{r^{2}}\right),
\]
 and in three dimensions,
\[
U_{ij}(\boldsymbol{r})=\frac{1}{8\pi}\left(-\frac{1}{r}\delta_{ij}+\frac{r_{i}r_{j}}{r^{3}}\right)
\]
where $\delta_{ij}$ the Kronecker delta. 

It can be shown that a generalized solution of the form 
\begin{equation}
\boldsymbol{u}(\boldsymbol{\alpha})=\int\boldsymbol{U}(\boldsymbol{x}(\boldsymbol{\alpha})-\boldsymbol{x}(\boldsymbol{\alpha}^{\prime}))\cdot\boldsymbol{f}(\boldsymbol{x}(\boldsymbol{\alpha}^{\prime}))d\boldsymbol{\alpha}^{\prime}\label{eq:GeneralizedSolution}
\end{equation}
is valid for $\boldsymbol{f}(\boldsymbol{x})$ in a wide variety of
function spaces \citep{ladyzhenskaya1969mathematical}.

In situations where $\boldsymbol{f}(\boldsymbol{x})$ is concentrated
on a lower dimensional manifold within a domain, Equation (\ref{eq:GeneralizedSolution})
remains valid \citep{cortez_method_2001}. For instance, if $\boldsymbol{f}(\boldsymbol{x})$
is concentrated on some closed surface, $S$, the integral equation
is equivalent to the single-layer hydrodynamic potential on the surface
$S$ of density $\boldsymbol{f}(\boldsymbol{x})$ \citep{ladyzhenskaya1969mathematical}. 

Of interest here and in many applications \citep{wrobel_enhanced_2016,wrobel_modeling_2014,aranda2015model},
is the quasi-steady state situation where the fluid velocity is found
from Stokes equation, but the boundary or particle positions and forces
may vary over time. This quasi-steady state assumption is applicable
for slowly moving immersed structures and small length scales. In
particular, the fluid velocity at any specific instant in time must
remain close to equilibrium velocity due to a boundary with the same
shape and force distribution as the moving boundary frozen at that
point in time. In such cases, the boundary position is typically updated
through the Eulerian-Lagrangian velocity equivalence: 
\[
\dot{\boldsymbol{x}}(\boldsymbol{\alpha},t)=\boldsymbol{u}(\boldsymbol{x}(\boldsymbol{\alpha},t),t)=\int\boldsymbol{U}(\boldsymbol{x}(\boldsymbol{\alpha},t)-\boldsymbol{x}(\boldsymbol{\alpha}^{\prime},t))\cdot\boldsymbol{f}(\boldsymbol{x}(\boldsymbol{\alpha}^{\prime},t),t)d\boldsymbol{\alpha}^{\prime}.
\]

The numerical discretization proceeds in two steps. The singular
fundamental solution is regularized by convolution with a smooth,
radially symmetric function that satisfies certain moment conditions
(e.g. \citep{cortez_method_2001,cortez_method_2005}), and the integral
is approximated as a summation over a finite number of points\footnote{For the quasi-steady state Stokes equations, the numerical discretization
procedure is quite similar to that used in the vortex method literature
\citep{anderson_vortex_1985,majda_vorticity_2002}.}, $\left\{ \boldsymbol{x}_{k}\right\} _{k=1}^{N}$, 
\begin{equation}
\dot{\boldsymbol{x}}_{k}(t)=\sum_{j=1}^{N}\boldsymbol{U}_{\epsilon}(\boldsymbol{x}_{k}(t)-\boldsymbol{x}_{j}(t))\boldsymbol{f}_{j}(\boldsymbol{x}_{j}(t),t)\Delta\alpha_{j}(t).\label{eq:RegularizedStokesletODE}
\end{equation}
The time dependence of $\Delta\alpha_{j}(t)$ occurs when the force
is concentrated on a lower dimensional manifold. In such cases, $\Delta\alpha_{j}(t)$
is treated as a surface area element which may vary in time as deformation
occurs. For the case of a collection of particles, $\Delta\alpha_{j}$
is generally a constant for each particle (e.g. the surface area,
or a drag-coefficient) such that $\boldsymbol{f}_{j}\Delta\alpha_{j}$
is equal to the force that the particle exerts on the fluid. Likewise,
for cases where the MRS is used to model fluid flow around a surface,
the regularization parameter $\epsilon$ typically should be dependent
on the mean surface area element size. On the other hand, for particle
simulations, $\epsilon$ typically corresponds to a physical length
scale related to the size of the particles \citep{cortez_method_2005}. 

Discussions of the well-posedness of Stokes equations and numerical
convergence results for the MRS can be found in \citep{cortez_method_2001,cortez_method_2005,ladyzhenskaya1969mathematical,kim_microhydrodynamics:_1991,pozrikidis1992boundary}.
 However, it appears that the well-posedness of the quasi-steady
state MRS for general forces has not been studied (though two recent
studies on specific examples in two dimensions exist \citep{mori_well-posedness_2017,lin_solvability_2017}).
On the other hand, for the spatially discrete systems considered in
this article, well-posedness is guaranteed by the Picard existence
and uniqueness theorem (for short times) as long as the force operator
is a Lipschitz continuous function of $\boldsymbol{x}(\boldsymbol{\alpha},t)$
and $t$ \citep[\S 3]{ladas1972differential}. 

\section{\emph{A Posteriori} Error Estimation\label{sec:A-Posteriori-Error}}

\emph{A posteriori} error estimation is a counterpart to the \emph{a
priori} error estimation techniques of classical numerical analysis.
With \emph{a priori} error estimation, bounds on the error in a numerical
method generally depend on derivatives of the exact solution to an
ODE or PDE and stability factors related to the numerical method.
A major difficulty with\emph{ a priori} bounds is that the exact solution
and its derivatives are typically unknown. This often prevents accurate
measures of the error in specific simulations from being obtained.
In contrast, \emph{a posteriori} error bounds do not depend upon
the analytical solution to a given problem. Instead, they typically
depend on finite differences of the numerical solution and stability
factors derived from the continuous ODE or PDE \citep{estep_posteriori_1995}.
We provide a brief overview of some of the results from \citep{estep_posteriori_1995}
regarding \emph{a posteriori }error estimation for ODEs. 

To start, we detail the finite element discretizations that will be
employed throughout the rest of the paper in Section \ref{subsec:Continuous-Galerkin-Finite}.
In Section \ref{subsec:Frchet-Derivatives}, we describe how the Fr\'echet
derivative of a nonlinear operator is obtained. This is needed to
form the linearized adjoint equations required by the \emph{a posteriori}
error estimation procedure to quantify how numerical errors propagate
over time. In Section \ref{subsec:Error-Representation-Formulas},
we derive error representation formulas for finite element methods
that give the numerical error in terms of computable quantities. In
Section \ref{subsec:Nodally-Equivalent-FEM}, we show how to derive
such error representation formulas for Runge-Kutta methods. Finally,
in Section \ref{subsec:ResidualQuadratureExplicit}, we discuss some
of the error components that we must monitor in the MRS. 

\subsection{Continuous Galerkin Finite Element Discretizations\label{subsec:Continuous-Galerkin-Finite}}

Consider the weak form of a system of ODEs of the form 
\begin{equation}
\begin{aligned}\left\langle \dot{\boldsymbol{x}}(t)-\boldsymbol{F}(t,\boldsymbol{x}(t)),\boldsymbol{V}(t)\right\rangle  & =0\,\,\,\,\forall\boldsymbol{V}\in\mathcal{V}\\
\left(\boldsymbol{x}(0),\boldsymbol{V}(0)\right) & =0
\end{aligned}
\label{eq:WeakForm}
\end{equation}
where $\mathcal{V}$ is a Hilbert space, and $\boldsymbol{F}(t,\boldsymbol{u})$
is assumed to be sufficiently smooth such that a solution, $\boldsymbol{x}(t)$,
to Equation (\ref{eq:WeakForm}) exists in some Hilbert space $\mathcal{U}.$
The symbol $\left\langle \cdot,\cdot\right\rangle $ is a duality
pairing on $\mathcal{V}$, and $\left(\cdot,\cdot\right)$ is an inner
product on $\mathbb{R}^{n}$. A Galerkin finite element method is
obtained by choosing finite dimensional spaces $\mathcal{V}_{n}\subset\mathcal{V}$
and $\mathcal{U}_{n}\subset\mathcal{U}$ such that there exists a
function $\boldsymbol{X}\in\mathcal{U}_{n}$ that satisfies 
\begin{align*}
\left\langle \dot{\boldsymbol{X}}(t)-\boldsymbol{F}(t,\boldsymbol{X}(t)),\boldsymbol{V}(t)\right\rangle  & =0\,\,\,\,\forall\boldsymbol{V}\in\mathcal{V}_{n}\\
\left(\boldsymbol{X}(0),\boldsymbol{V}(0)\right) & =0
\end{align*}
Functions in the space $\mathcal{V}_{n}$ are called \emph{test functions},
and functions in $\mathcal{U}_{n}$ are called \emph{trial functions}.
Typically, for finite element methods, the domain, $[0,T]$ is partitioned
into a finite number of subintervals, $[0,T]=\cup_{i=1}^{N}\mathcal{I}_{i}$
where $\mathcal{I}_{i}=[t_{i-1},t_{i})$. Refinements of this partition
then involve the addition of partition points $t_{i}$ into the domain.
For a continuous Galerkin method of order $q$, we then construct
$\mathcal{U}_{N}$ and $\mathcal{V}_{N}$ as spaces of piecewise continuous
polynomials. In particular, we define $\mathcal{P}^{q}(\mathcal{I}_{i};\mathbb{R}^{n})$
as the set of polynomials of order less than or equal to $q$ with
domain $\mathcal{I}_{i}$ and range $\mathbb{R}^{n}$. Then we have:
\[
\mathcal{U}_{N}=\left\{ \boldsymbol{u}(t)|\forall i\in[1,N],\,\,\boldsymbol{u}|_{\mathcal{I}_{i}}\in\mathcal{P}^{q}(\mathcal{I}_{i};\mathbb{R}^{n})\mbox{ and }\boldsymbol{u}\in C^{0}([0,T];\mathbb{R}^{n})\right\} 
\]
\[
\mathcal{V}_{N}=\left\{ \boldsymbol{V}(t)|\forall i\in[1,N],\,\,\boldsymbol{V}|_{\mathcal{I}_{i}}\in\mathcal{P}^{q-1}(\mathcal{I}_{i};\mathbb{R}^{n})\right\} .
\]
 where $\boldsymbol{u}|_{\mathcal{I}}$ and $\boldsymbol{V}|_{\mathcal{I}}$
are the restrictions of $\boldsymbol{u}$ and $\boldsymbol{V}$ to
an interval $\mathcal{I}\subseteq[0,T]$.

\subsection{Fr\'echet Derivatives of Operators\label{subsec:Frchet-Derivatives}}

Let $F:Y\rightarrow Z$ be an operator on $Y\subset X$ where $X$
and $Z$ are Banach spaces and $Y$ is an open subset of $X$. This
operator is Fr\'echet differentiable if there exists a bounded linear
operator $DF$ such that 
\[
F[u+\delta u]-F[u]=(DF[u])\delta u+o(\Vert\delta u\Vert).
\]

For any particular choice of $u(x)\in Y$ the form of this operator
can often be found explicitly by considering the limit 
\[
DF[u]v=\lim_{\epsilon\rightarrow0}\frac{F[u+\epsilon v]-F[u]}{\epsilon\Vert v\Vert}
\]
for some arbitrary $v\in Y$. 

\subsection{Derivation of Error Representation Formulas\label{subsec:Error-Representation-Formulas}}

Follwing \citep{estep_posteriori_1995}, the starting point for both
\emph{a priori} and \emph{a posteriori} analysis is to subtract the
exact and discrete equations and linearize about the range of trajectories
from $\boldsymbol{x}(t)$ to $\boldsymbol{X}(t)$. We define the error
$\boldsymbol{e}(t)=\boldsymbol{x}(t)-\boldsymbol{X}(t)$ and introduce
the operator, 
\[
B(\boldsymbol{W},\boldsymbol{V})\equiv\left\langle \dot{\boldsymbol{W}}(t)-\int_{0}^{1}D\boldsymbol{F}(t,s\boldsymbol{x}(t)+(1-s)\boldsymbol{X}(t))ds\cdot\boldsymbol{W}(t),\boldsymbol{V}(t)\right\rangle 
\]
where $D\boldsymbol{F}(\cdot)$ is the Fr\'echet derivative of $\boldsymbol{F}(t,\boldsymbol{x})$
with respect to variations in $\boldsymbol{x}$.  For short-hand,
we denote the linearized operator, $D\boldsymbol{F}(\cdot)$ averaged
over all trajectories between the continuum and numerical solutions
as 
\[
\boldsymbol{A}(t)=\int_{0}^{1}D\boldsymbol{F}(t,s\boldsymbol{x}(t)+(1-s)\boldsymbol{X}(t))ds,
\]
and integrating by parts, we obtain
\[
B(\boldsymbol{W},\boldsymbol{V})=\left\langle \boldsymbol{W}(t),-\dot{\boldsymbol{V}}(t)-\boldsymbol{A}^{\ast}(t)\boldsymbol{V}(t)\right\rangle -(\boldsymbol{W}(t_{n}),\boldsymbol{V}(t_{n}))\,\,\,\,\,\,\forall\boldsymbol{V}\in\mathcal{V}^{q-1}([0,T]).
\]
Choosing $\boldsymbol{Z}(t)$ that solves
\[
B(\boldsymbol{W},\boldsymbol{Z})=(\boldsymbol{W}(T),\boldsymbol{Z}(T))\,\,\,\,\,\forall\boldsymbol{W}\in\mathcal{V}^{q}([0,T]),
\]
and substituting $\boldsymbol{e}$ for $\boldsymbol{W},$ we obtain
the \emph{a priori} error representation formula, 
\[
B(\boldsymbol{e},\boldsymbol{Z})=|\boldsymbol{e}(T)|^{2}.
\]
Using Galerkin orthogonality, and the linearity of $B(\cdot,\boldsymbol{Z})$,
\begin{equation}
B(\boldsymbol{e}-\pi\boldsymbol{e},\boldsymbol{Z})=|\boldsymbol{e}(T)|^{2}\label{eq:a_priori_rep}
\end{equation}
where $\pi\boldsymbol{e}$ is the $L^{2}$ projection of the error
into $\mathcal{U}_{N}$. Hence, $\boldsymbol{Z}(t)$ is the continuous
Galerkin approximation in $\mathcal{V}^{q}([0,T])$ of the weak solution
of the continuous adjoint problem
\begin{equation}
\begin{aligned}-\boldsymbol{\dot{z}}(t)-\boldsymbol{A}^{\ast}(t)\boldsymbol{z}(t) & =0\\
\boldsymbol{z}(T) & =\boldsymbol{e}(T)/\left|\boldsymbol{e}(T)\right|.
\end{aligned}
\label{eq:AdjointContinuous}
\end{equation}
Further analysis of Equation (\ref{eq:a_priori_rep}) leads to bounds
on the error that depend on derivatives of $\boldsymbol{x}(t)$ and
the stability properties of the numerical scheme (through $\boldsymbol{Z}(t)$).
On the other hand, if we solve Equation (\ref{eq:AdjointContinuous})
for $\boldsymbol{z}(t)$, apply Galerkin orthogonality, and compute
\[
\left\langle \dot{\boldsymbol{X}}(t)-\boldsymbol{F}(t,\boldsymbol{X}(t)),\boldsymbol{z}(t)-\pi\boldsymbol{z}(t)\right\rangle =|\boldsymbol{e}(t_{n})|^{2}
\]
we obtain bounds that depend on $\boldsymbol{X}(t)$ and stability
factors which are functionals of $\boldsymbol{z}(t)$ \citep{estep_posteriori_1995}. 

For finite element method discretizations, where the numerical solution
is defined at every point, \emph{a posteriori} analysis can be conducted
in a well-defined way using weak formulations over continuous and
discrete function spaces. On the other hand, for finite difference
methods, where the solution is only defined on a finite set of points,
more work must be done to define a suitable adjoint equation and residual
operator. One way of approaching this issue is through the construction
of special finite element methods that are related to the finite difference
method in question. This is described in the next subsection.

\subsection{Nodally Equivalent Finite Element Methods\label{subsec:Nodally-Equivalent-FEM}}

Because Runge-Kutta methods do not fall into the scope of typical
finite element analysis, the first step to obtaining an error representation
formula is the development of a \emph{nodally equivalent finite element
method }(neFEM)\citep{collins_posteriori_2015} that allows us to
consistently extrapolate the pointwise Runge-Kutta solution (defined
at time steps $t_{n}$) to a globally defined, piecewise polynomial
solution defined on $[0,T].$ This is accomplished by starting with
a standard continuous Galerkin method, and then applying various projection
operators and quadrature rules so that at the time nodes, $t_{n}$,
the solution of the modified continuous Galerkin method is equivalent
to that of the Runge Kutta method. 

In general, an $L$ stage Runge Kutta method can be written in the
form
\begin{equation}
\begin{array}{clll}
\boldsymbol{k}_{\ell} & = & \boldsymbol{F}\left(t_{n}+c_{\ell}\Delta t_{n},\boldsymbol{X}_{n}+h_{n}\sum_{j=0}^{L}a_{\ell j}\boldsymbol{k}_{j}\right) & \forall\ell\in[1,2,\dots,L]\\
\boldsymbol{X}_{n+1} & = & \boldsymbol{X}_{n}+\Delta t\sum b_{\ell}\boldsymbol{k}_{\ell}
\end{array}\label{eq:RungeKuttaMethod}
\end{equation}
For an explicit Runge-Kutta method, $a_{\ell j}=0$ whenever $j\geq\ell$.
In this paper, only explicit Runge-Kutta methods are considered. 

The construction of neFEM\emph{ }is discussed in \citep{collins_posteriori_2015}.
To summarize, a neFEM is formed by introducing certain projection
operators and numerical quadratures to the weak formulation of the
problem.The projection operators, $\mathcal{P}_{\ell}$, correspond
to the extrapolations done at each stage of an explicit Runge-Kutta
method and can be written in the form 
\[
\mathcal{P}_{\ell}[\boldsymbol{X}](t)=\boldsymbol{X}_{n}(t_{n})+\boldsymbol{k}_{\ell}(t-t_{n})
\]
where $\boldsymbol{k}_{\ell}$ is the $\ell$th stage of the Runge-Kutta
method. With this formulation, $\boldsymbol{k}_{\ell}$ can be written
in terms of $\boldsymbol{F}(t,\mathcal{P}_{k}\boldsymbol{X})$ with
$k<\ell$ in the explicit Runge-Kutta case. This leads to a modified
variational formulation, 
\[
\left\langle \dot{\boldsymbol{X}}-\sum_{\ell}\boldsymbol{F}\left(t,\mathcal{P}_{\ell}\left[\boldsymbol{X}\right](t)\right),\boldsymbol{V}\right\rangle =0\,\,\,\,\,\forall\boldsymbol{V}\in\mathcal{V}^{q-1}([0,T]).
\]
Next, quadrature rules are used to approximate integrals of the form
\[
\int_{t_{n}}^{t_{n+1}}\boldsymbol{F}(t,\mathcal{P}_{\ell}[\boldsymbol{X}](t))dt.
\]
As an example, we may use the midpoint rule, 
\[
\int_{t_{n}}^{t_{n+1}}\boldsymbol{F}(t,\mathcal{P}_{\ell}[\boldsymbol{X}](t))dt\approx\Delta t\boldsymbol{F}\left[\frac{1}{2}\left(t_{n}+t_{n+1}\right),\mathcal{P}_{\ell}[\boldsymbol{X}]\left(\frac{1}{2}\left(t_{n}+t_{n+1}\right)\right)\right].
\]
Combining the multiple stages of a Runge-Kutta method together, we
write
\[
\left\langle \dot{\boldsymbol{X}},\boldsymbol{V}\right\rangle =\sum_{\ell}\left\langle \boldsymbol{F}(t,\mathcal{P}_{\ell}\boldsymbol{X}),\boldsymbol{V}\right\rangle _{d_{\ell}}\,\,\,\,\forall\boldsymbol{V}\in\mathcal{V}^{q-1}(I_{n})
\]
where $\left\langle \cdot,\cdot\right\rangle $ represents an integral,
and $\left\langle \cdot,\cdot\right\rangle _{d_{\ell}}$ a numerical
quadrature rule. Choosing a basis (e.g. Legendre polynomials) for
$\mathcal{V}^{q-1}(I_{n})$ allows us to solve for some polynomial
$\boldsymbol{X}(t)=\sum_{k=0}^{q-1}\boldsymbol{a}_{k}(t-t_{n})^{k}$
over each interval $[t_{n},t_{n+1}]$ which defines the solution at
all time points. 

\subsection{Residual, Quadrature, and Explicit Errors\label{subsec:ResidualQuadratureExplicit}}

In \citep{collins_posteriori_2015}, \emph{a posteriori} estimation
techniques are applied to explicit multistep and Runge-Kutta time
stepping methods. With implicit finite element methods, the error
terms that result from \emph{a posteriori} analysis involve residuals,
\begin{equation}
\boldsymbol{R}[\boldsymbol{X}]\equiv\dot{\boldsymbol{X}}-\boldsymbol{F}(t,\boldsymbol{X}).\label{eq:Residual}
\end{equation}
 However, when explicit methods are used, there are two additional
sources of error. In our formulation of Runge-Kutta methods as neFEM
methods, discrete quadrature approximations introduce quadrature errors,
and the approximation of $\boldsymbol{X}(t)$ by various projections
introduces an extrapolation error. We briefly summarize Theorem 3
from \citep{collins_posteriori_2015} which shows why these error
terms appear. 

With the introduction of the projection operators described in the
previous section, the ``continuous'' residual of Equation (\ref{eq:Residual})
is modified as in \citep{collins_posteriori_2015} with $b_{\ell}$
as in Equation (\ref{eq:RungeKuttaMethod})
\[
\boldsymbol{R}_{P}[\boldsymbol{X}]=\dot{\boldsymbol{X}}-\sum_{\ell}b_{\ell}\boldsymbol{F}(t,\mathcal{P}_{\ell}\boldsymbol{X}).
\]

Along with the introduction of numerical quadrature and a choice of
finite dimensional function space, this leads to a modified version
of the standard Galerkin orthogonality,
\[
\left\langle \dot{\boldsymbol{X}},\boldsymbol{V}\right\rangle -\sum_{\ell=1}^{L}b_{\ell}\left\langle \boldsymbol{F}(t,\mathcal{P}_{\ell}\boldsymbol{X}),\boldsymbol{V}\right\rangle _{d_{\ell}}=0\,\,\,\,\,\forall\boldsymbol{V}\in\mathcal{V}^{q-1}(I_{i}),\,\,i\in[1,2,\dots,m].
\]
Setting $\boldsymbol{V}(t)=\boldsymbol{z}(t)$ where $\boldsymbol{z}(t)$
is the solution to Equation (\ref{eq:AdjointContinuous}), and adding
and subtracting\\  $\sum_{\ell=1}^{L}b_{\ell}\boldsymbol{F}(t,\mathcal{P}_{\ell}\boldsymbol{X})$,
we obtain 
\[
\left\langle \boldsymbol{R}[\boldsymbol{X}],\boldsymbol{z}\right\rangle =\left\langle \boldsymbol{R}_{P}[\boldsymbol{X}],\boldsymbol{z}\right\rangle +\left\langle \sum_{\ell}^{L}b_{\ell}\boldsymbol{F}(t,\mathcal{P}_{\ell}\boldsymbol{X})-\boldsymbol{F}(t,\boldsymbol{X}),\boldsymbol{z}\right\rangle 
\]
Now, applying the modified Galerkin orthogonality, 
\[
\left\langle \boldsymbol{R}[\boldsymbol{X}],\boldsymbol{z}\right\rangle =\boldsymbol{E}_{R}[\boldsymbol{X},\boldsymbol{z}]+\boldsymbol{E}_{E}[\boldsymbol{X},\boldsymbol{z}]+\boldsymbol{E}_{Q}\left[\boldsymbol{X},\boldsymbol{z}\right]
\]
with 
\begin{align}
\boldsymbol{E}_{R}[\boldsymbol{X},\boldsymbol{z}] & \equiv\left\langle \boldsymbol{R}_{P}[\boldsymbol{X}],\boldsymbol{z}-\pi\boldsymbol{z}\right\rangle \label{eq:ResidualError}\\
\boldsymbol{E}_{E}[\boldsymbol{X},\boldsymbol{z}] & \equiv\left\langle \sum_{\ell}^{L}b_{\ell}\boldsymbol{F}(t,\mathcal{P}_{\ell}\boldsymbol{X})-\boldsymbol{F}(t,\boldsymbol{X}),\boldsymbol{z}\right\rangle \label{eq:ExplicitError}\\
\boldsymbol{E}_{Q}\left[\boldsymbol{X},\boldsymbol{z}\right] & \equiv\sum_{\ell=1}^{L}b_{\ell}\left(\left\langle \boldsymbol{F}(t,\mathcal{P}_{\ell}\boldsymbol{X}),\pi\boldsymbol{z}\right\rangle _{d_{\ell}}-\left\langle \boldsymbol{F}(t,\mathcal{P}_{\ell}\boldsymbol{X}),\pi\boldsymbol{z}\right\rangle \right)\label{eq:QuadratureError}
\end{align}
 Equation (\ref{eq:ResidualError}) is a residual error that measures
how well the ODE can be approximated in the finite dimensional space,
$\mathcal{V}^{q}([0,T])$, Equation (\ref{eq:ExplicitError}) is an
explicit error term resulting from the extrapolation of $\boldsymbol{X}$
across each interval, and Equation (\ref{eq:QuadratureError}) is
a quadrature error term. A key point in this last step is that the
operator $(\mathcal{I}-\pi)$ can often be bounded in terms of derivatives.
For instance, bounds of the form $\left|\boldsymbol{z}-\pi\boldsymbol{z}\right|\leq Ch\left|\dot{\boldsymbol{z}}\right|$
are frequently found \citep{estep_posteriori_1995}. 

\section{Application to the Method of Regularized Stokeslets }

Recalling that $\Delta\alpha_{j}$ is either a surface area element,
or a conserved physical quantity, such as the mass of a particle in
a fluid, we introduce the shorthand, 

\[
\mathcal{S}_{\epsilon}[\boldsymbol{x}]_{k}\equiv\sum_{j}\boldsymbol{U}_{\epsilon}(\boldsymbol{x}_{k}-\boldsymbol{x}_{j})\boldsymbol{F}_{j}[\boldsymbol{x}]\Delta\boldsymbol{\alpha}_{j}
\]
so that the method of regularized Stokeslets equations can be written
as $\dot{\boldsymbol{x}}_{k}=\mathcal{S}_{\epsilon}[\boldsymbol{x}]_{k}$.

\subsection{The Adjoint of the Spatially Discrete MRS Operator}

For the MRS applied to a network of particles the Fr\'echet derivative
is of the form
\begin{align*}
D\boldsymbol{\mathcal{S}}_{\epsilon}^{h}[\boldsymbol{x}](\boldsymbol{y})= \sum_{j}\left(\nabla_{x}\boldsymbol{U}_{\epsilon}\right)(\boldsymbol{x}_{k}-\boldsymbol{x}_{j}):\left(\boldsymbol{F}_{j}[\boldsymbol{x}]\otimes(\boldsymbol{y}_{k}-\boldsymbol{y}_{j})\right)h^{2}+\sum_{j}\boldsymbol{U}_{\epsilon}(\boldsymbol{x}_{k}-\boldsymbol{x}_{j})\cdot\nabla_{x}\boldsymbol{F}_{j}[\boldsymbol{x}]\cdot\boldsymbol{y}_{j}h^{2}.
\end{align*}
where $\otimes$ is a dyadic product\footnote{The dyadic product is defined as $\boldsymbol{a}\otimes\boldsymbol{b}=\boldsymbol{a}\boldsymbol{b}^{T}$
with $\boldsymbol{a}$ and $\boldsymbol{b}$ in $\mathbb{R}^{3}$ }, and $:$ denotes a second order tensor contraction. The adjoint,
can be obtained by considering the inner product 
\[
\left\langle D\boldsymbol{\mathcal{S}}_{\epsilon}^{h}[\boldsymbol{x}](\boldsymbol{y}),\boldsymbol{\phi}\right\rangle =\sum_{k}\sum_{j}\left(\nabla_{x}\boldsymbol{U}_{\epsilon}\right)(\boldsymbol{x}_{k}-\boldsymbol{x}_{j}):\left(\boldsymbol{F}_{j}[\boldsymbol{x}]\otimes(\boldsymbol{y}_{k}-\boldsymbol{y}_{j})\right)h^{2}+\sum_{j}\boldsymbol{U}_{\epsilon}(\boldsymbol{x}_{k}-\boldsymbol{x}_{j})\cdot\nabla_{x}\boldsymbol{F}_{j}[\boldsymbol{x}]\cdot\boldsymbol{y}_{j}h^{2}\boldsymbol{\phi}_{k}
\]
where $\boldsymbol{\phi}$ is contained in the same function space
as $\boldsymbol{y}$. To obtain an adjoint operator, we must isolate
all terms that depend on $\boldsymbol{y}_{k}$ for some particular
$k$. We note that in this case, $\nabla_{x}\boldsymbol{K}(\boldsymbol{x}_{k}-\boldsymbol{x}_{j})=-\nabla_{x}\boldsymbol{K}(\boldsymbol{x}_{j}-\boldsymbol{x}_{k})$
and that the product adjoint formula: $(AB)^{T}=B^{T}A^{T}$ should
be applied when necessary. The adjoint operator is of the form
\begin{align*}
&\left(D\boldsymbol{\mathcal{S}}_{\epsilon}^{h}[\boldsymbol{x}]\right)_{k}^{\ast}(\boldsymbol{\phi})=\\ &\sum_{j}\left(\nabla_{x}\boldsymbol{U}_{\epsilon}\right)^{T}(\boldsymbol{x}_{k}-\boldsymbol{x}_{j}):(\boldsymbol{F}_{j}[\boldsymbol{x}]\otimes\boldsymbol{\phi}_{k}+\boldsymbol{\phi}_{j}\otimes\boldsymbol{F}_{k}[\boldsymbol{x}])h^{2}+\nabla_{x}\boldsymbol{F}_{k}[\boldsymbol{x}]^{T}:\left(\sum_{j}\boldsymbol{U}_{\epsilon}(\boldsymbol{x}_{k}-\boldsymbol{x}_{j})\cdot\boldsymbol{\phi}_{j}h^{2}\right).
\end{align*}
As an example, if each pair of attached points are connected by an
elastic spring, we can write
\[
\boldsymbol{F}_{j}[\boldsymbol{x}]=\sum_{\mathcal{N}_{i}}k\left(\frac{\left|\boldsymbol{x}_{i}-\boldsymbol{x}_{j}\right|}{r_{ij}^{0}}-1\right)\frac{\boldsymbol{x}_{i}-\boldsymbol{x}_{j}}{\left|\boldsymbol{x}_{i}-\boldsymbol{x}_{j}\right|},
\]
 where the summation is over points connected to $\boldsymbol{x}_{j}$,
and $r_{ij}^{0}$ is the equilibirium separation of $\boldsymbol{x}_{i}$
and $\boldsymbol{x}_{j}$. The Fr\'echet derivative of this operator
is 
\[
D\boldsymbol{F}[\boldsymbol{x}](\boldsymbol{y})=\sum_{\mathcal{N}_{i}}k\left[\left(\frac{1}{r_{ij}^{0}}-\frac{1}{\left|\boldsymbol{x}_{i}-\boldsymbol{x}_{j}\right|}\right)\left(\boldsymbol{\mathcal{I}}-\frac{\boldsymbol{x}_{i}-\boldsymbol{x}_{j}}{\left|\boldsymbol{x}_{i}-\boldsymbol{x}_{j}\right|}\otimes\frac{\boldsymbol{x}_{i}-\boldsymbol{x}_{j}}{\left|\boldsymbol{x}_{i}-\boldsymbol{x}_{j}\right|}\right)+\frac{1}{r_{ij}^{0}}\frac{\boldsymbol{x}_{i}-\boldsymbol{x}_{j}}{\left|\boldsymbol{x}_{i}-\boldsymbol{x}_{j}\right|}\otimes\frac{\boldsymbol{x}_{i}-\boldsymbol{x}_{j}}{\left|\boldsymbol{x}_{i}-\boldsymbol{x}_{j}\right|}\right]\cdot(\boldsymbol{y}_{i}-\boldsymbol{y}_{j}).
\]
If the connectivity matrix describing the connections between $\boldsymbol{x}_{j}$
and $\boldsymbol{x}_{k}$ is symmetric (as it typically will be due
to physical considerations), then this operator is self-adjoint. 

\subsection{Regularization Error}

In addition to the residual, quadrature, and explicit errors, there
is a regularization error term associated with the MRS. This term
arises from the use of regularized Stokeslets instead of singular
Stokeslets in Equation (\ref{eq:RegularizedStokesletODE}). It can
be written as 
\begin{equation}
\boldsymbol{E}_{Re}\left[\boldsymbol{X},\boldsymbol{z}\right]\equiv\left\langle \boldsymbol{\mathcal{S}}_{\epsilon}^{h}[\boldsymbol{X}]-\boldsymbol{\mathcal{S}}_{0}^{h}[\boldsymbol{X}],\boldsymbol{z}\right\rangle .\label{eq:RegularizationError}
\end{equation}
We also note that consideration of regularization error effects how
the adjoint is defined. If the continuous problem involves no regularization,
then the adjoint is defined with respect to a modified problem (due
to the use of Runge-Kutta methods), but without regularization as
\begin{equation}
\dot{\boldsymbol{x}}-\sum_{\ell}b_{\ell}\boldsymbol{\mathcal{S}}_{0}^{h}[\mathcal{P}_{\ell}\boldsymbol{x}].\label{eq:ModifiedProblem}
\end{equation}
In \citep[\S 5]{collins_posteriori_2015}, the effect of instabilities
caused by numerical approximation of an operator on the choice of
adjoint is discussed. In this case, since regularization leads to
a more stable system, it does not introduce instabilities that are
artifacts of the numerics. Thus, the complicated ``dual-adjoint''
procedure discussed in \citep[\S 5]{collins_posteriori_2015} is not
needed here. 

Furthermore, although the operator, $\boldsymbol{\mathcal{S}}_{0}^{h}[\mathcal{P}_{\ell}\boldsymbol{x}]$
is an unbounded operator on $\boldsymbol{x}$, once $\boldsymbol{x}$
has been found, the linearized adjoint of $\boldsymbol{\mathcal{S}}_{0}^{h}[\mathcal{P}_{\ell}\boldsymbol{x}]$
may be treated as a bounded linear operator acting on the adjoint
solution. 

\subsection{Numerical Error Estimation and Adjoint Equation Solution Algorithm}

With the error representation formula and adjoint equation of the
previous section, we now discuss how the error terms are computed
once a numerical approximation is found and the adjoint equation solved.
All three error terms, Equations (\ref{eq:ResidualError})-(\ref{eq:QuadratureError})
of Section \ref{subsec:ResidualQuadratureExplicit} depend upon integration
over time of a complicated function. These integrals are not analytical
in general, and are approximated by numerical quadrature, or by deriving
a bound on the size of the operator. 

Even with quadrature approximations no reference to the continuous
solution is needed to obtain bounds on the error in the quadrature
used to approximate the error terms. The use of quadrature risks introducing
unreliability into the error bounds, since quadrature may underestimate
the integral\footnote{Quadrature may also cause overestimation of the integral. In this
case, error bound remains valid, but loses its sharpness}. However, in most situations, the use of high order quadrature is
likely to be accurate so long as the numerical solution is sufficiently
smooth, and not severely under-resolved. Furthermore, since the error
in quadrature approximations can be bounded in terms of derivatives
of the integrand, which in this case is a function of the numerical
solution, this approach preserves the \emph{a posteriori} nature of
the bounds. 

In practice, given some explicit Runge-Kutta method of order $p$,
we use Gaussian quadrature formulas to approximate the time integration.
Since the projection operators of the \emph{neFEM} are defined for
any time, it is possible to evaluate $\boldsymbol{\mathcal{S}}_{\epsilon}^{h}\left[\mathcal{P}_{\ell}\boldsymbol{x}\right]$
at any time, thus any standard quadrature method may be employed.
We use sufficiently high order Gaussian quadrature formulas so that
the error in estimating the integrals is much smaller than the value
of the integral. 

Now that we have discussed how the error terms can be estimated,
it remains to discuss how to solve the adjoint equation for $\boldsymbol{z}(t)$.
For the forward problems, we use Runge-Kutta methods of order 4 or
less. In the \emph{a posteriori }it is customary to use a higher order
solver for the adjoint equation in order to ensure that the error
terms are accurately computed. . Thus, we use a 6th order Runge-Kutta
method \citep{butcher2009fifth}. In many practical applications,
$\boldsymbol{z}$ is only needed for controlling the error and is
not as important in terms of physical insight as the numerical approximant,
$\boldsymbol{x}$. Thus, it may be the case that crude estimates for
$\boldsymbol{z}$ are sufficient. 

The ``correct'' initial condition (or final condition since data
is specified at $t=T$) is to set $\boldsymbol{\boldsymbol{z}}(T)=\boldsymbol{e}(T)/\Vert\boldsymbol{e}(T)\Vert$.
This is problematic since $\boldsymbol{e}(T)$ is unknown. However,
approximations of the initial condition for the adjoint equation were
studied in \citep{cao_posteriori_2004}. It was shown that randomized
initial conditions for the adjoint equation perform well for error
estimation, and that taking the maximum error approximation over several
initial values will lead, with increasingly high probability, to nearly
optimal results. Since we lack information about what $\boldsymbol{z}(T)$
should be, we chose $\boldsymbol{z}(T)$ to be a random vector of
unit norm. For the case where the ODE is a spatial discretization
of a PDE, we attempt to more accutely capture spatial correlations
in the error by drawing initial data from a Gaussian spatial process.
Such initial data can be obtained efficiently using Fourier transform
techniques discussed in \citep{kroese2013spatial}.

\section{Numerical Results\label{sec:Numerical-Results}}

In computing the forward solution, we use several common Runge-Kutta
methods. In particular, we use the RK4 method and the second order
accurate Heun's method. We track over each interval, the values of
$\boldsymbol{x}(t_{n})$ and also $\boldsymbol{k}_{i}$ $i=1,\dots\ell$,
the result from each stage of the RK method. These values are sufficient
to reconstruct a polynomial FEM solution over each interval in a partition
of $[0,T]$. For the adjoint equation we obtain a solution through
the use of a sixth order Runge-Kutta method \citep{butcher2009fifth}.
It is important to note that although the Runge Kutta method may exhibit
a certain order of convergence at the time nodes, the extrapolated
polynomial solution may be of lower degree over the interval. For
instance with the RK4 method the extrapolated polynomial solution
is accurate to second order. This fact is verified through our numerical
results; when the time step is halved, the error is quartered, but
superconvergence occurs at the time nodes due to special cancellations
of error \citep{collins_posteriori_2015}. As is done in \emph{a priori
}numerical estimation, to ensure that the Runge Kutta methods obtain
the theoretically expected convergence rates, we compute convergence
factors of the form, 
\[
\rho=\frac{\log\vert\boldsymbol{X}_{h}(T)-\tilde{\boldsymbol{x}}(T)\vert}{\log\vert\boldsymbol{X}_{h/2}(T)-\tilde{\boldsymbol{x}}(T)\vert}
\]
where $\tilde{\boldsymbol{x}}$ is a highly refined numerical approximation
to $\boldsymbol{x}$, and $\boldsymbol{X}_{h}$ is the numerical approximation
with time step $h$. For the RK4 method, it was observed that $\rho\approx4$
and for the RK6 method $\rho\approx6$ for the problems we tested. 

\subsection{Method of Regularized Stokeslets Example }

In this case we consider a deformed circle as the initial starting
position and consider the errors developed for a tethered boundary
and an elastic boundary. We set the resting configuration of the boundary
to be a circle of radius 1, and the initial position of the boundary
is parametrized as
\[
\boldsymbol{x}_{0}(s)=\left(\cos(\pi s)+\frac{1}{2}\sin\left(2\pi\cos\left(\pi(s-1)\right)\right),\,\sin(\pi s)\right)\,\,\,\,\,s\in[0,1)
\]
and simulate the boundary motion as it deforms towards its resting
position. 

\begin{figure}
\begin{centering}
\subfloat[]{\begin{centering}
\includegraphics[width=0.3\columnwidth]{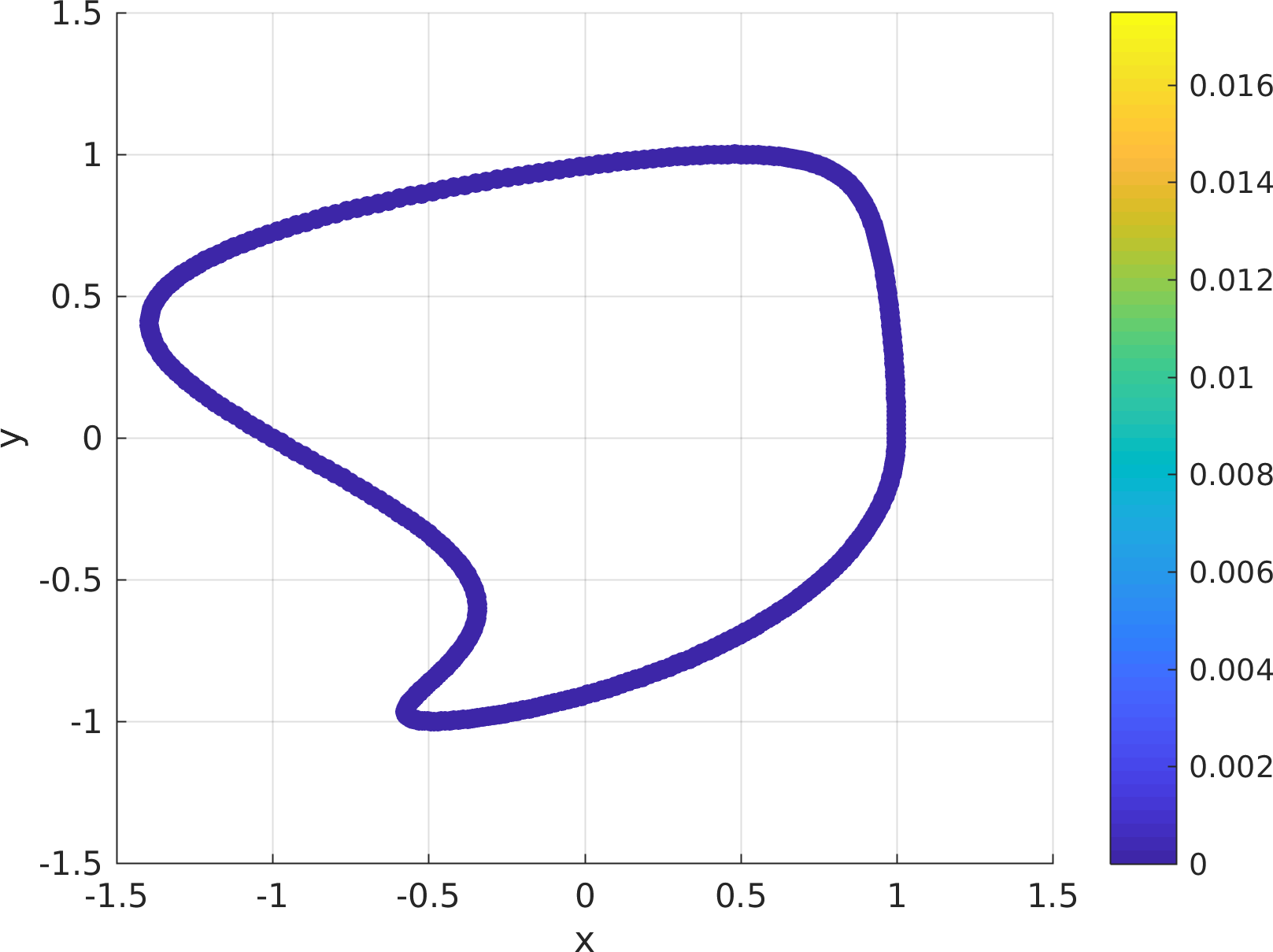}
\par\end{centering}
}\enskip{}\subfloat[]{\begin{centering}
\includegraphics[width=0.3\columnwidth]{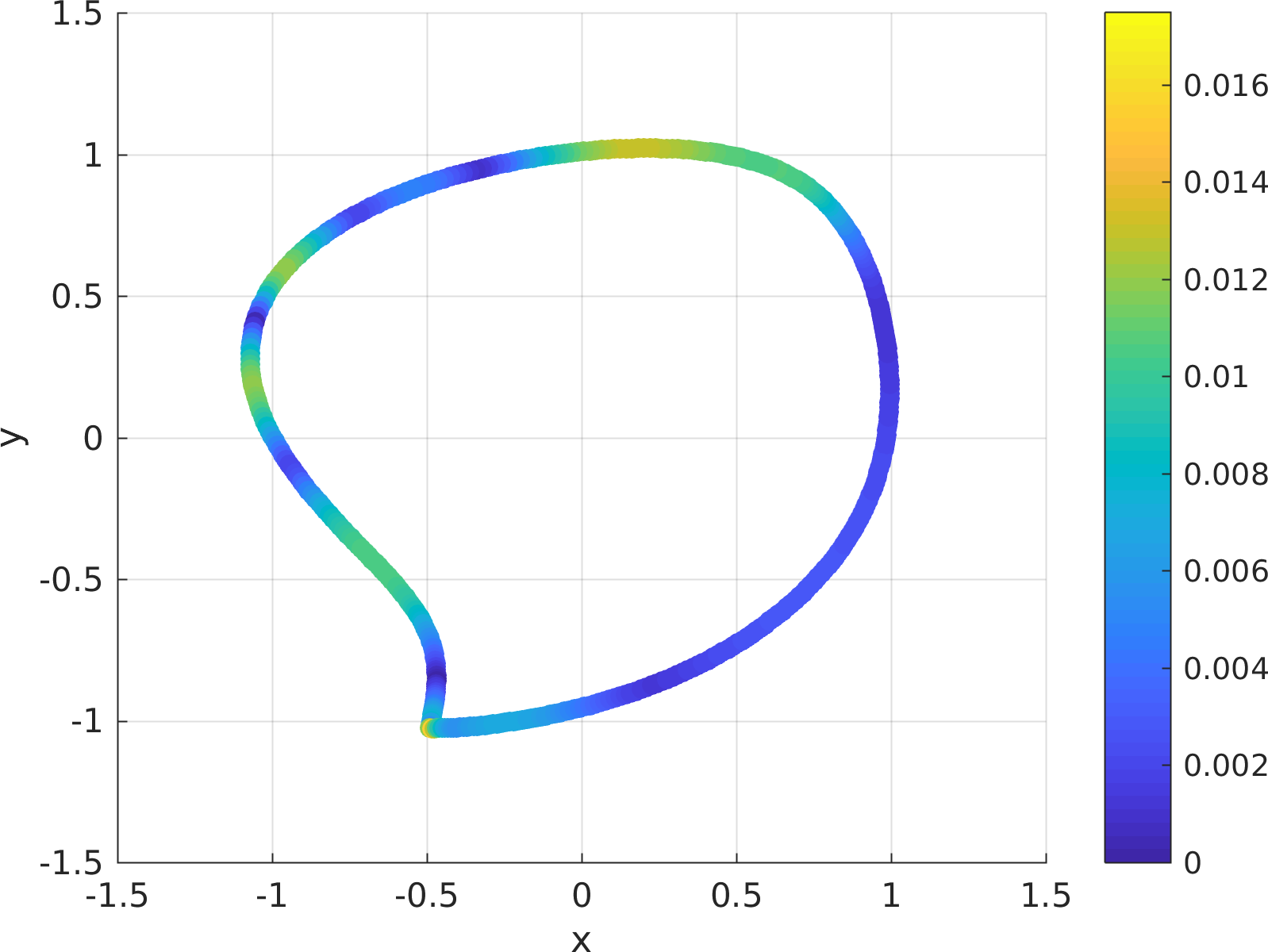}
\par\end{centering}
}\enskip{}\subfloat[]{\begin{centering}
\includegraphics[width=0.3\columnwidth]{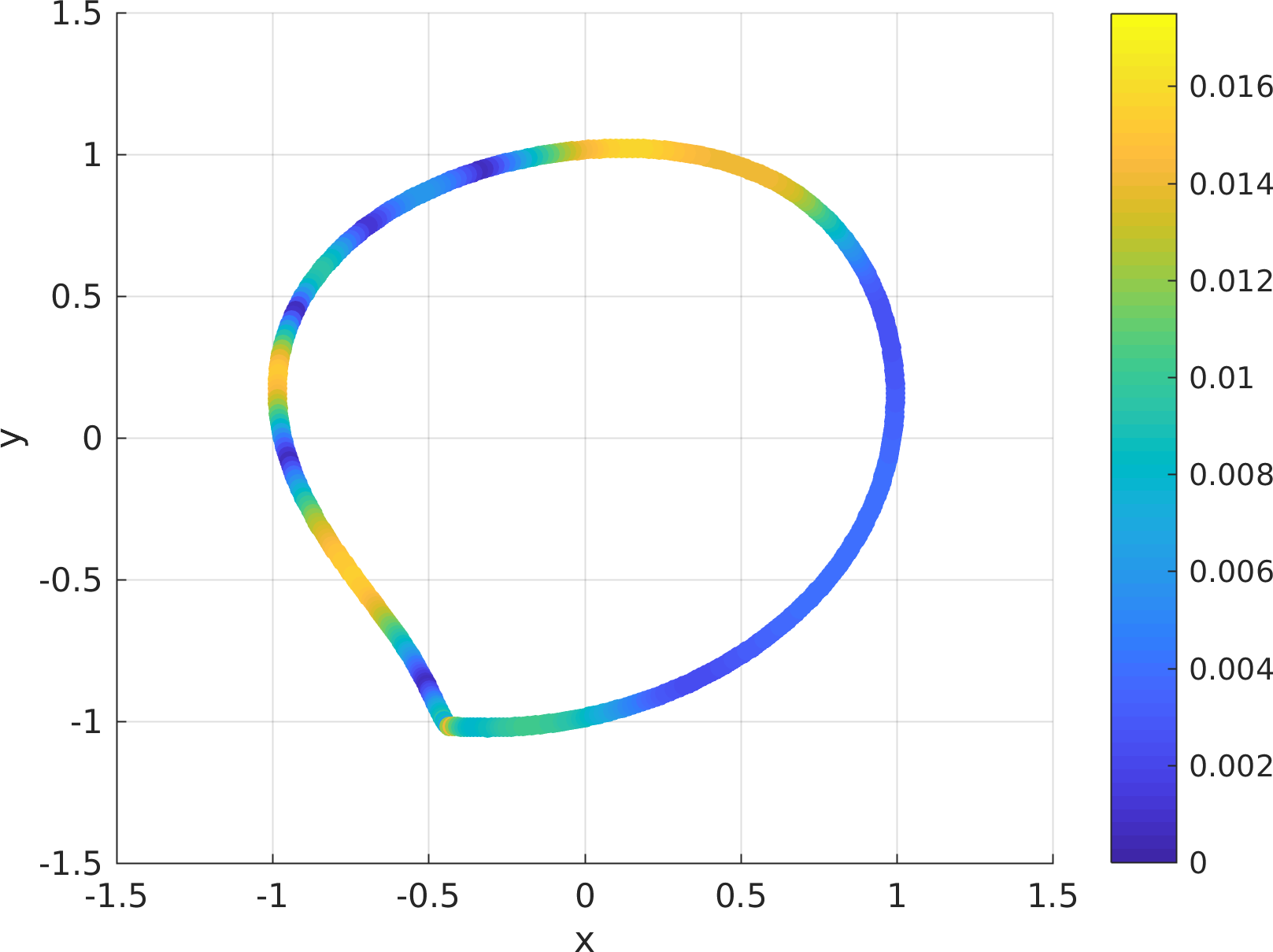}
\par\end{centering}
}
\par\end{centering}
\caption{\label{fig:StokesletsError}The regularization error term and position
of boundary shown at $t=0.0s$ (a), $t=0.5s$ (b), and $t=1.0s$ (c). }
\end{figure}

\begin{figure}
\begin{centering}
\subfloat[]{\begin{centering}
\includegraphics[width=0.3\columnwidth]{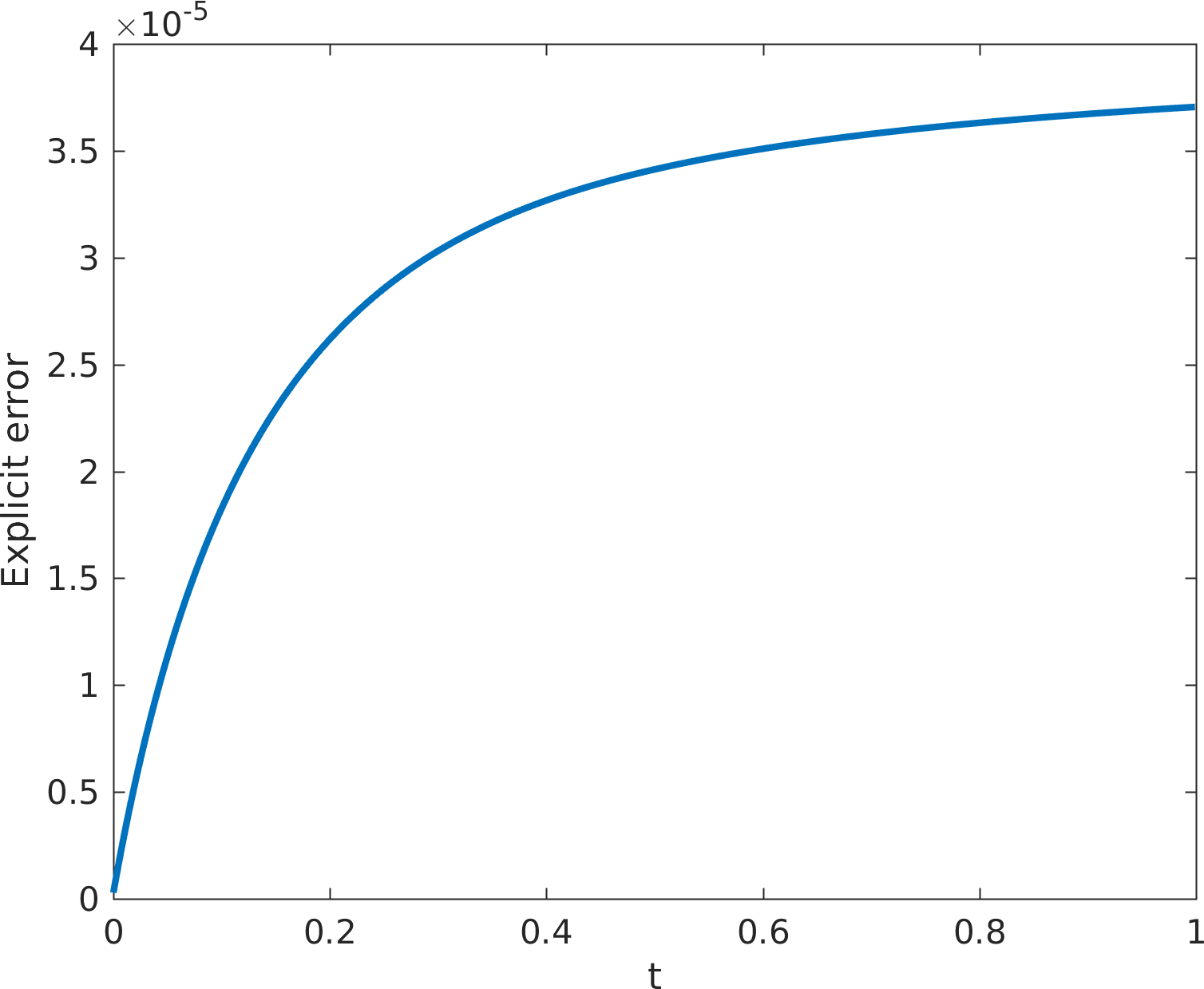}
\par\end{centering}
}\enskip{}\subfloat[]{\begin{centering}
\includegraphics[width=0.3\columnwidth]{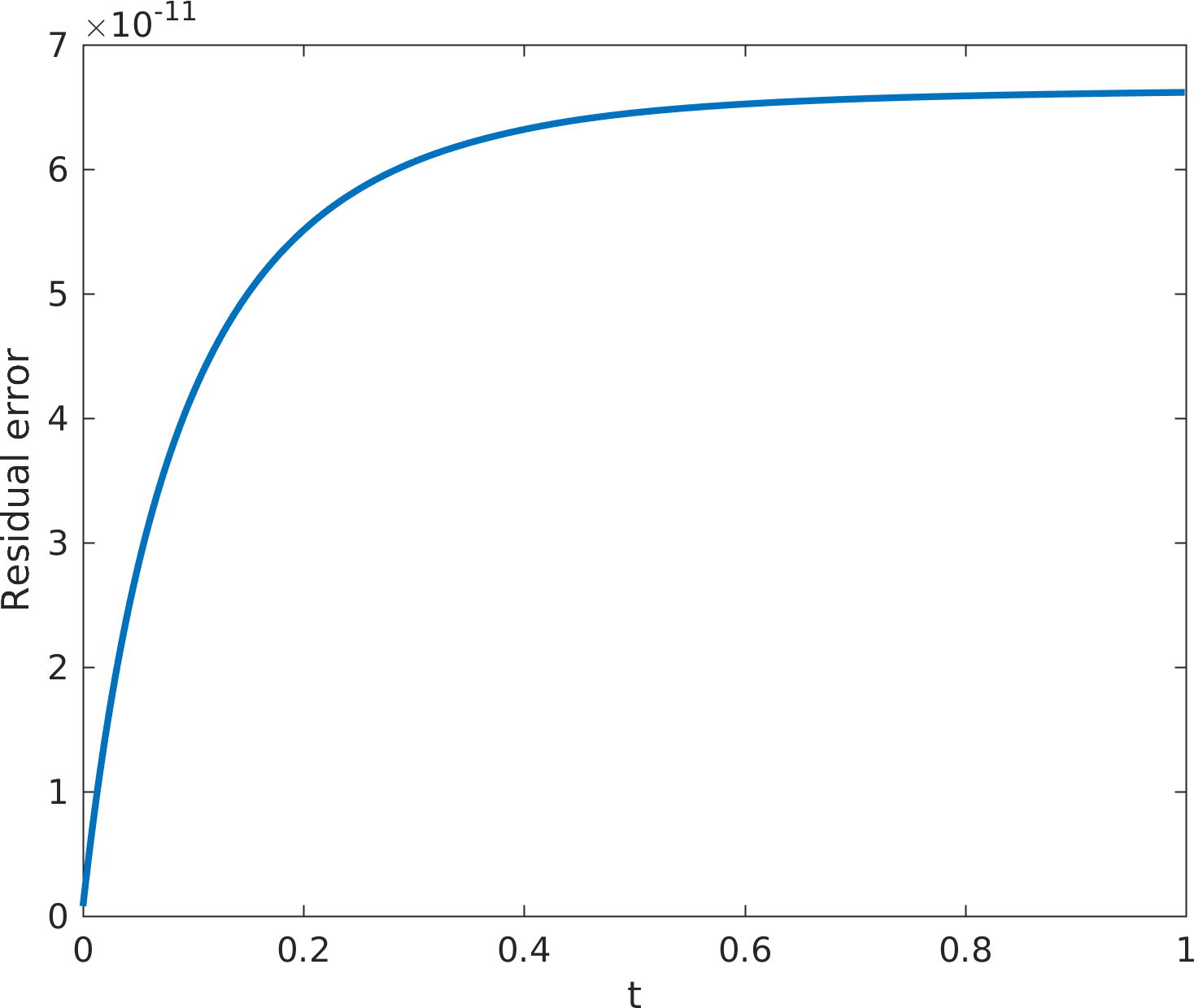}
\par\end{centering}
}\\
\subfloat[]{\begin{centering}
\includegraphics[width=0.3\columnwidth]{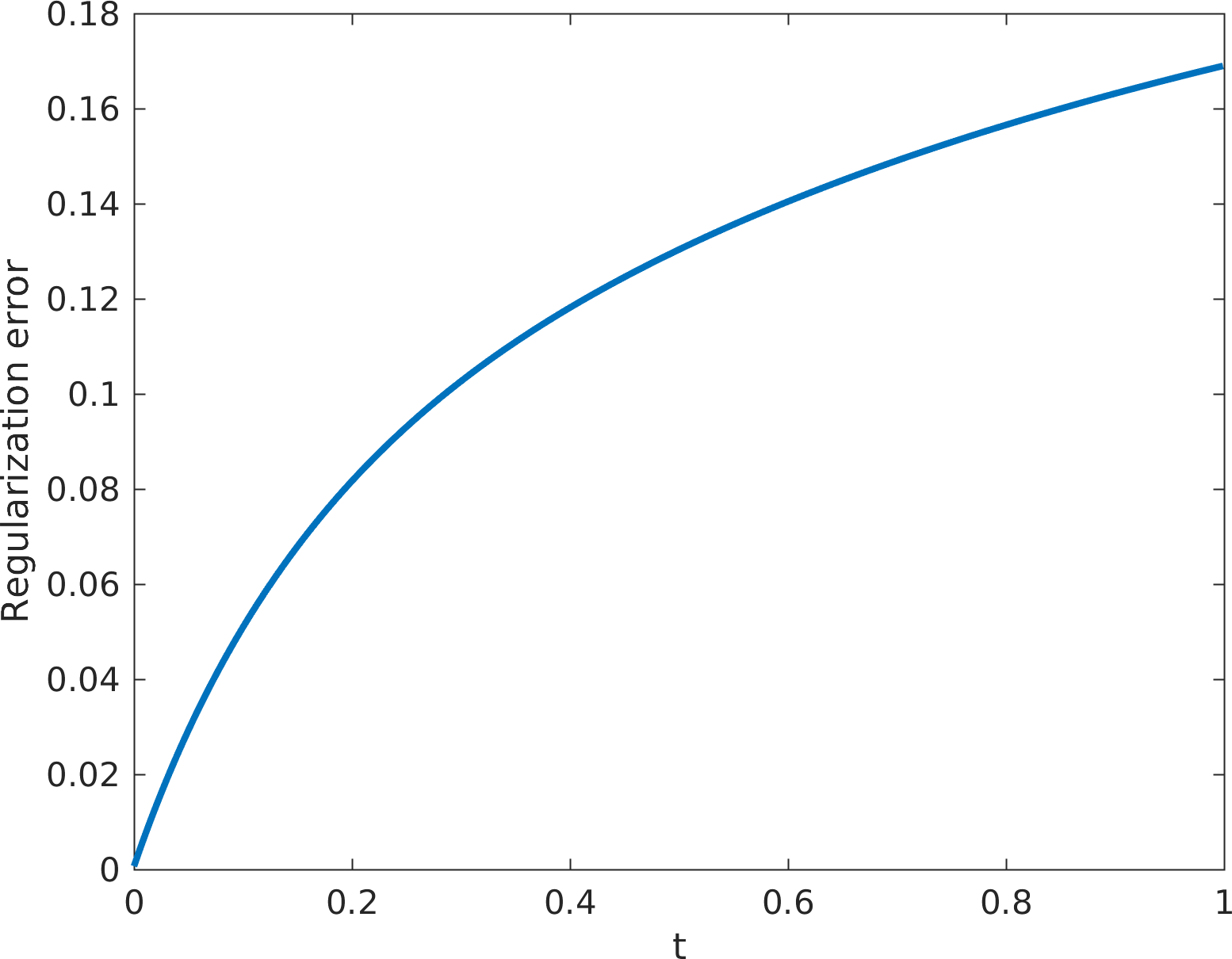}
\par\end{centering}
}\enskip{}\subfloat[]{\begin{centering}
\includegraphics[width=0.3\columnwidth]{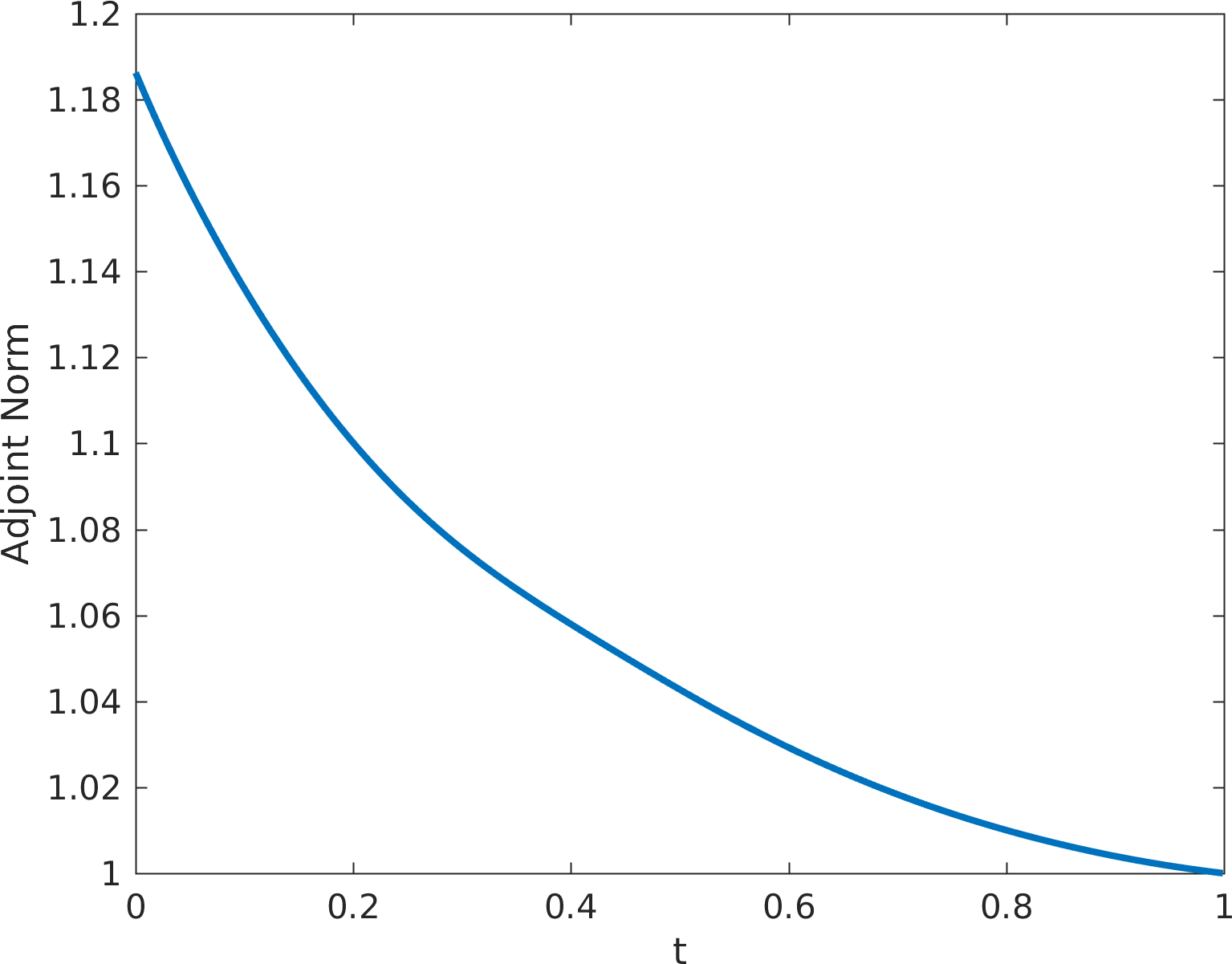}
\par\end{centering}
}
\par\end{centering}
\caption{\label{fig:ErrorPropagation}Explicit Error, Residual Error, Regularization
Error, and the norm of $\boldsymbol{z}(\cdot,t)$ as functions of
time. }
\end{figure}

In Figure \ref{fig:StokesletsError}, we see the residual error and
adjoint solution to the problem described above. We see that the initially
deformed shape relaxes towards its resting configuration as a circle.
We also observe some concentration in the error where the curvature
is highest. In Figure \ref{fig:ErrorPropagation}, the time dependence
of the explicit error, residual error, and adjoint solution norms
are shown. The quadrature error is not shown since it behaves very
similarly to the explicit error term. The error terms grow rapidly
at first, followed by a period of slow growth as the boundary approaches
its resting configuration. The adjoint solution remains small even
at $t=0$ and gradually decays across the time interval.

\begin{figure}
\begin{centering}
\subfloat[]{\begin{centering}
\includegraphics[width=0.3\columnwidth]{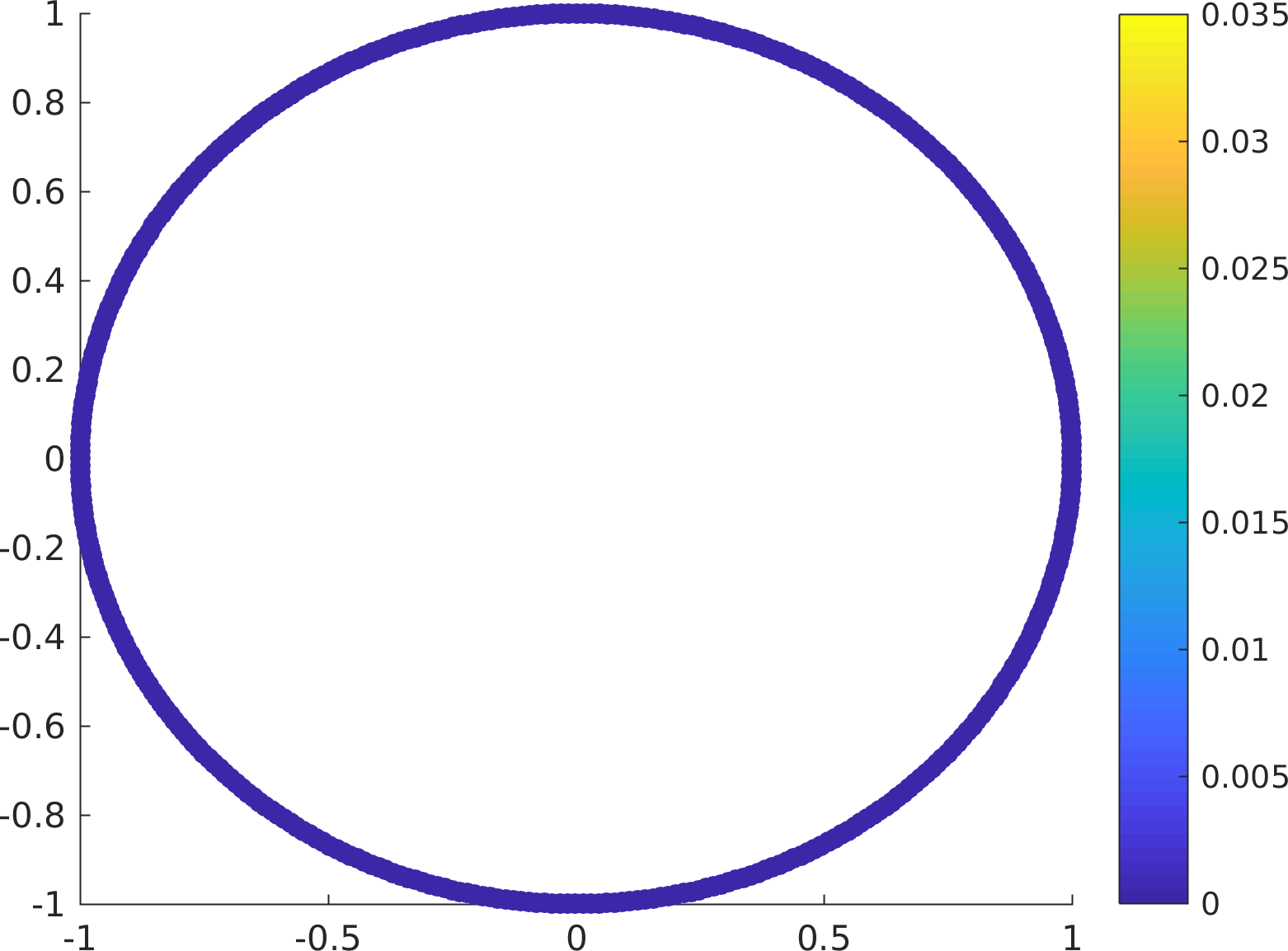}
\par\end{centering}
}\enskip{}\subfloat[]{\begin{centering}
\includegraphics[width=0.3\columnwidth]{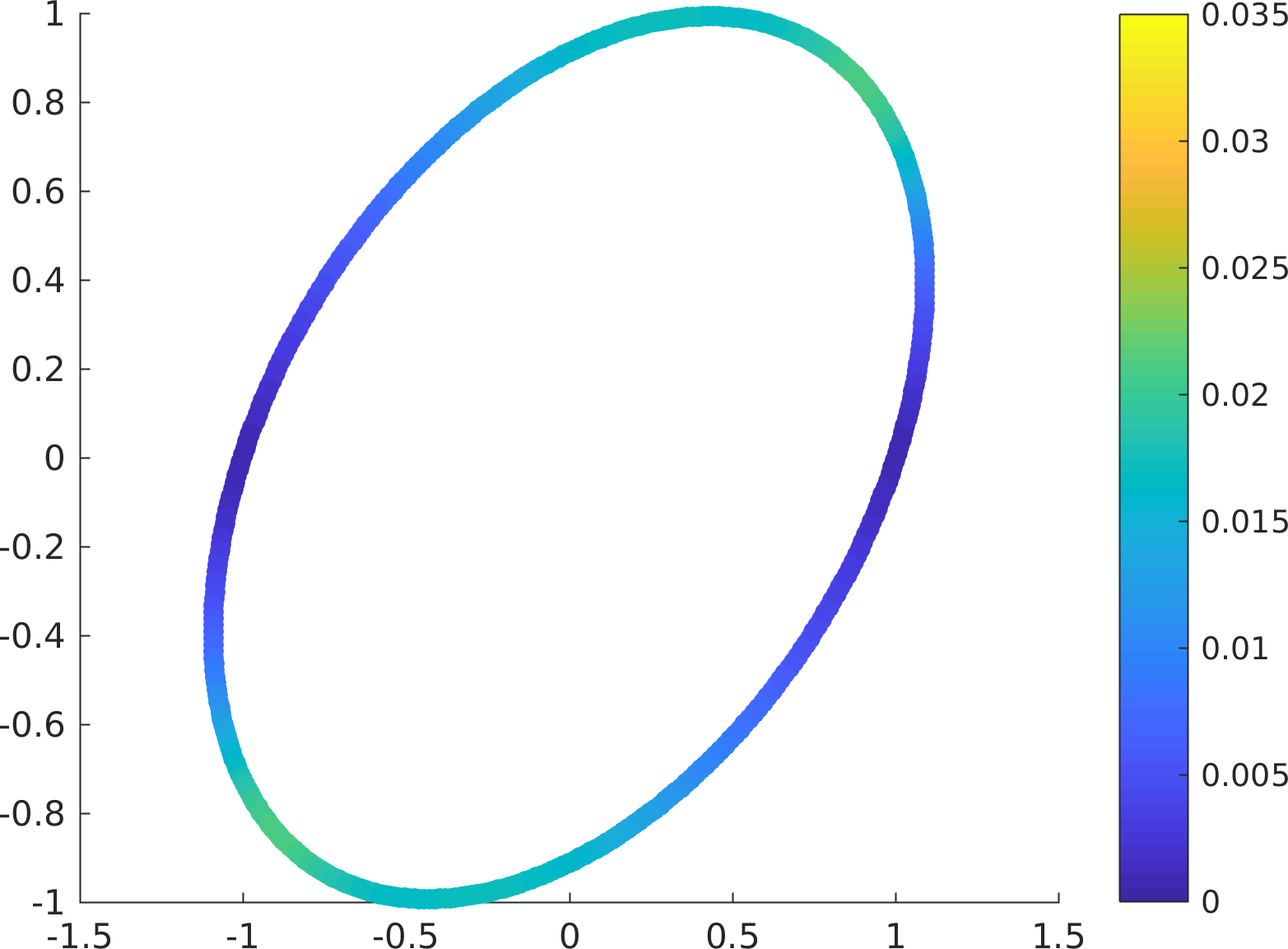}
\par\end{centering}
}\enskip{}\subfloat[]{\begin{centering}
\includegraphics[width=0.3\columnwidth]{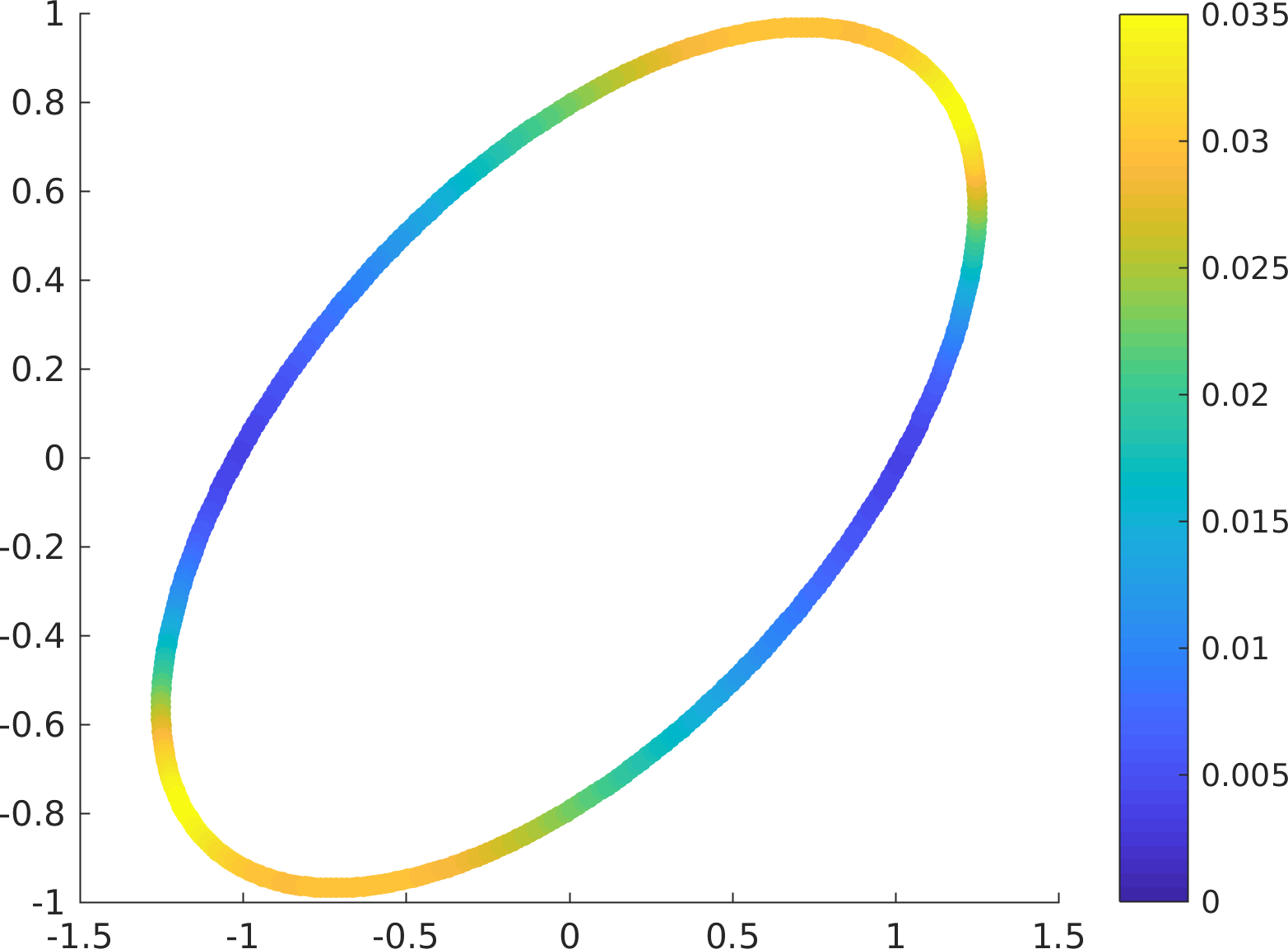}
\par\end{centering}
}
\par\end{centering}
\caption{\label{fig:StokesletsError2}The regularization error term and position
of boundary shown at $t=0.0s$ (a), $t=0.5s$ (b), and $t=1.0s$ (c).}
\end{figure}

\begin{figure}
\begin{centering}
\subfloat[]{\begin{centering}
\includegraphics[width=0.3\columnwidth]{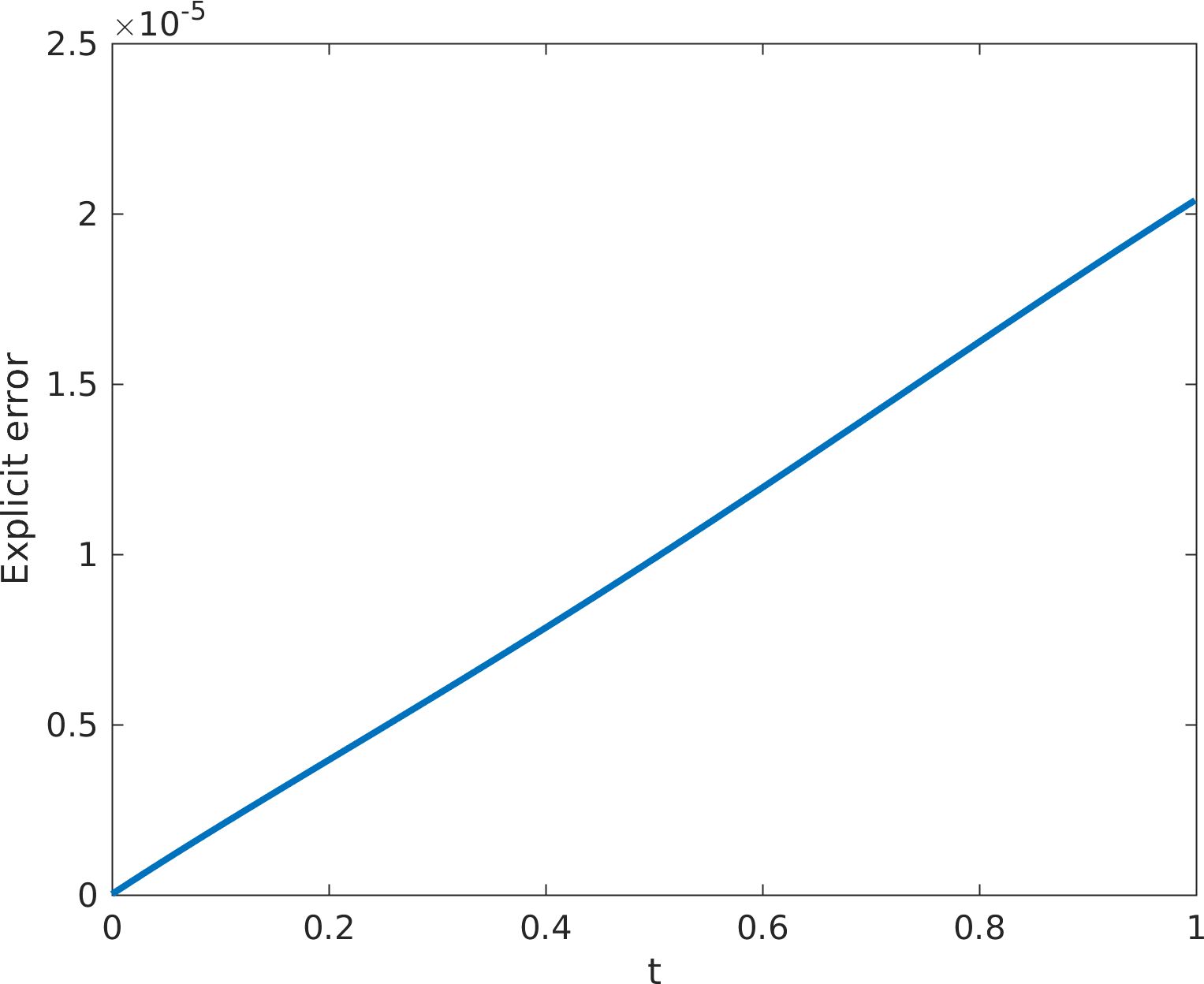}
\par\end{centering}
}\enskip{}\subfloat[]{\begin{centering}
\includegraphics[width=0.3\columnwidth]{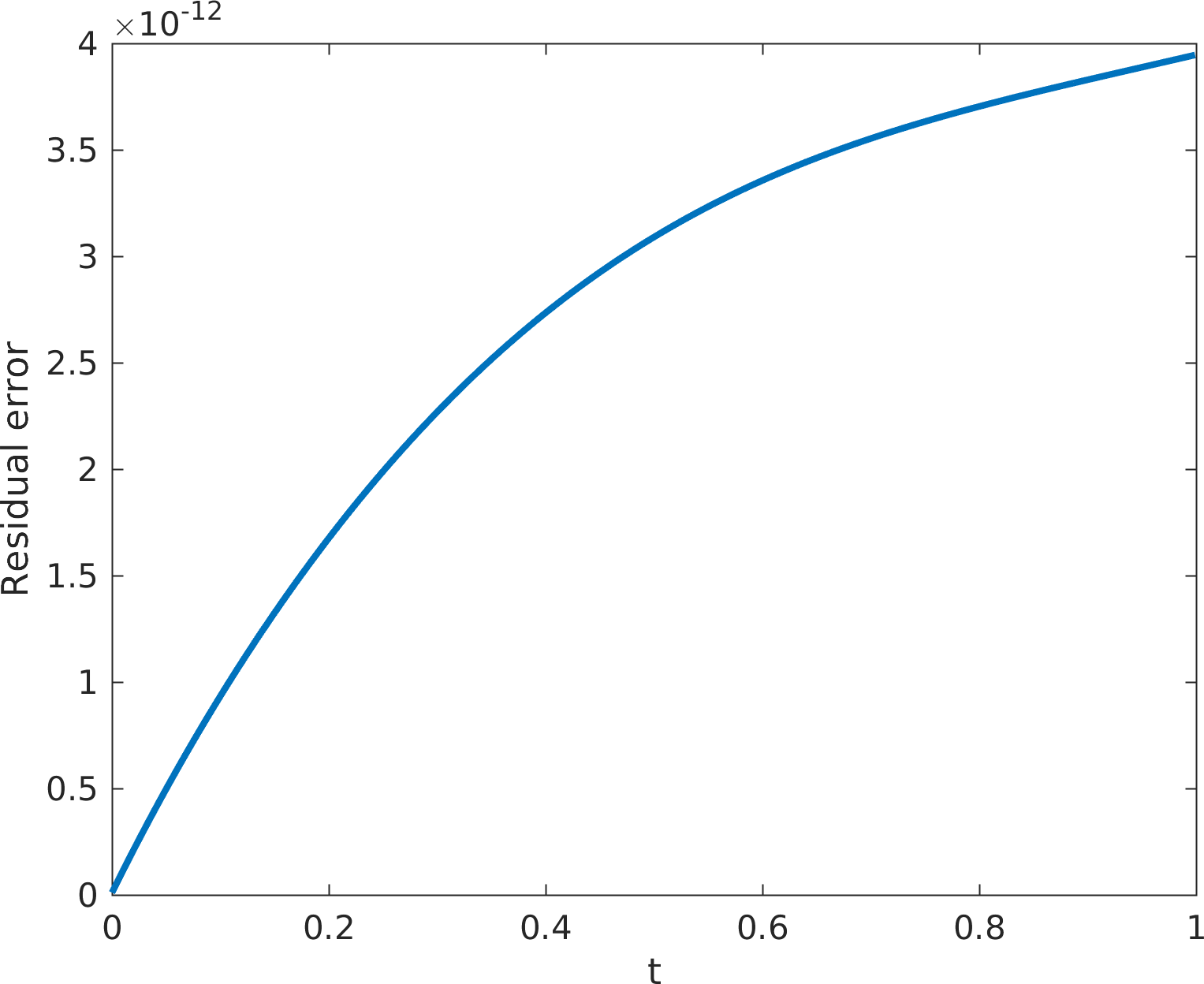}
\par\end{centering}
}\\
\subfloat[]{\begin{centering}
\includegraphics[width=0.3\columnwidth]{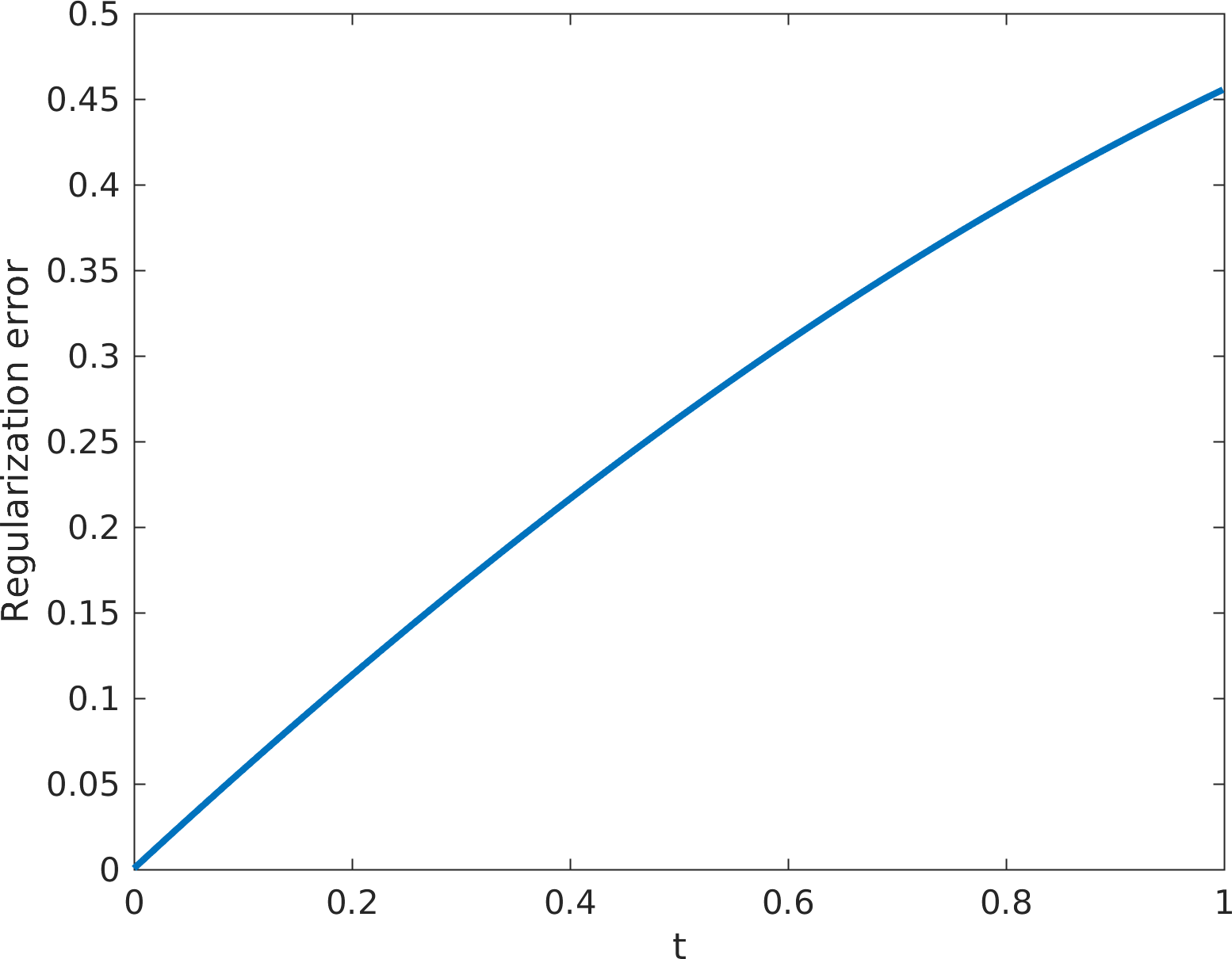}
\par\end{centering}
}\enskip{}\subfloat[]{\begin{centering}
\includegraphics[width=0.3\columnwidth]{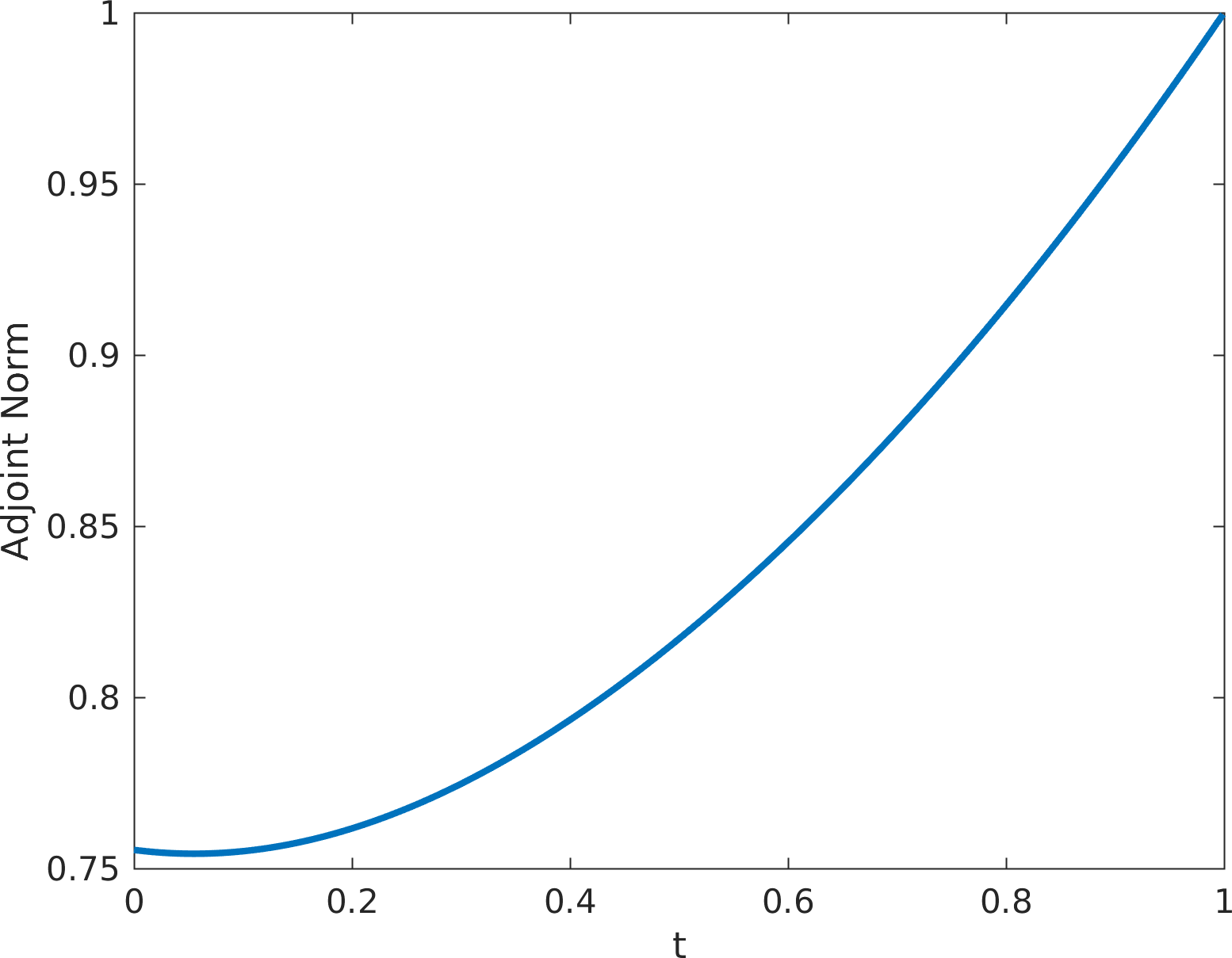}
\par\end{centering}
}
\par\end{centering}
\caption{\label{fig:ErrorPropagation2}Explicit Error, Residual Error, Regularization
Error, and the norm of $\boldsymbol{z}(\cdot,t)$ as functions of
time. }
\end{figure}

As a second example, we consider a circular boundary that is deformed
in a background shear velocity field of the form $\boldsymbol{u}=\left(y,0\right)$.
In Figure \ref{fig:StokesletsError2}, we see the extension and rotation
of a circular elastic body in response to a shear flow field. The
coloring of the deforming object corresponds to the magnitude of the
regularization error. Figure \ref{fig:ErrorPropagation2} depicts
the time dependence of the $\ell_{2}$ norm over space during this
deformation. We see that the adjoint solution remains small, but that
given the numerical parameters of our simulation, the regularization
error dominates the residual and explicit errors and increases roughly
linearly in time. 


\subsection{Method of Regularized Stokeslets applied to an Elastic Network of
Fibers}

In \citep{wrobel_modeling_2014}, the method of regularized Stokeslets
was used to model a three dimensional network of fibers immersed in
a fluid. This is of interest as there are number of application in
biology where the fluid structure interactions of such materials are
of importance. The application in \citep{wrobel_modeling_2014} was
geared towards understanding how a spermatocyte swims through a material
known as the zona pellucida that surrounds ooctyes. Another potential
application is in the study of biofilms growing in a slowly moving
fluid. In this case, the particles of the method of regularized Stokeslets
represent bacteria cells, and are not elements of a discretized surface.
In this case, the ODE system is exact no longer an approximation of
a PDE, but forms the governing equations. Regularization in this
context is typically used as a means of ensuring stability of the
resulting numerical simulations by limiting the velocity when bacteria
approach each other. The choice of length scale then depends on the
level of resolution needed in the simulation. 

Following \citep{wrobel_modeling_2014,wrobel_enhanced_2016}, we write
\begin{eqnarray*}
\dot{\boldsymbol{x}}_{k}(t) & = & \sum_{j}\boldsymbol{U}_{\epsilon}(\boldsymbol{x}_{k}-\boldsymbol{x}_{j})\boldsymbol{f}_{kj}(t)m_{j}\\
\boldsymbol{f}_{kj}(t) & = & \ell_{kj}^{2}E_{kj}\left(\frac{\vert\boldsymbol{x}_{k}(t)-\boldsymbol{x}_{j}(t)\vert}{\ell_{kj}}-1\right)\frac{\boldsymbol{x}_{k}(t)-\boldsymbol{x}_{j}(t)}{\vert\boldsymbol{x}_{k}(t)-\boldsymbol{x}_{j}(t)\vert}
\end{eqnarray*}
Note that $f_{kj}(t)=-f_{jk}(t)$. In \citep{wrobel_enhanced_2016,wrobel_modeling_2014},
the resting length is allowed to change in order to incorporate viscoelastic
effects, but for simplicity, we consider only the elastic (constant
resting length) case here. 

\begin{figure}
\begin{centering}
\subfloat[]{\begin{centering}
\includegraphics[width=0.3\columnwidth]{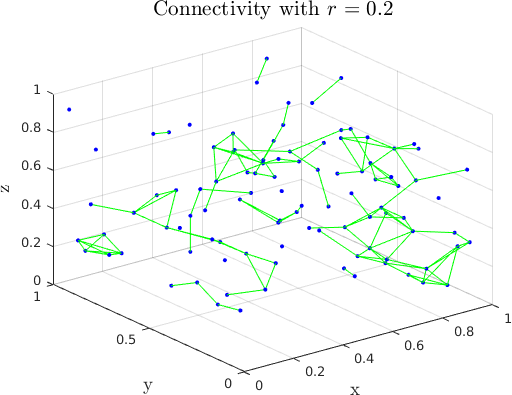}
\par\end{centering}
}\enskip{}\subfloat[]{\begin{centering}
\includegraphics[width=0.3\columnwidth]{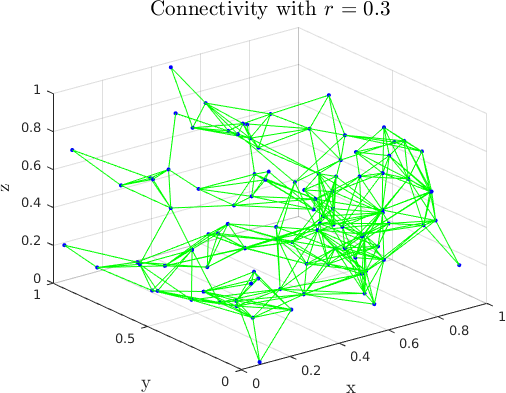}
\par\end{centering}
}\enskip{}\subfloat[]{\begin{centering}
\includegraphics[width=0.3\columnwidth]{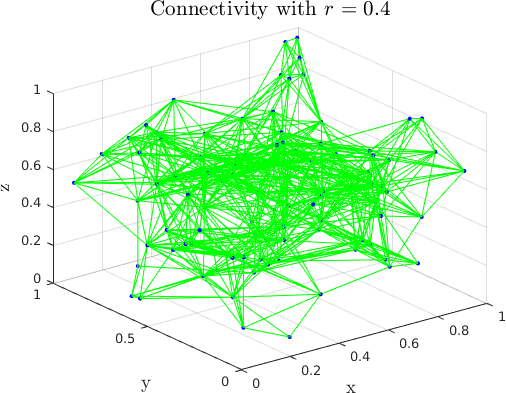}
\par\end{centering}
}
\par\end{centering}
\caption{\label{fig:A-network-of-fibers}A network of elastic fibers with 100
points for different values for $r_{connect}$. }
\end{figure}

As a test problem, we consider 100 points that form a network of fibers
with a connectivity rule that if $\vert\boldsymbol{x}_{i}-\boldsymbol{x}_{j}\vert\leq r_{connect}$,
at time 0, then $\boldsymbol{x}_{i}$ and $\boldsymbol{x}_{j}$ are
attached by an elastic spring, and if $\vert\boldsymbol{x}_{i}-\boldsymbol{x}_{j}\vert>r_{connect}$,
then there is no attachment. In Figure \ref{fig:A-network-of-fibers},
an image of such a network is depicted. 

As an example problem, we place the elastically connected bacteria
in a linear flow field, $\boldsymbol{u}(\boldsymbol{x})=\boldsymbol{C}\cdot\boldsymbol{x}$
where $C$ is a second order tensor. The velocity of each point is
then set to 
\[
\dot{\boldsymbol{x}}_{k}(t)=\sum_{j}\boldsymbol{U}_{\delta}(\boldsymbol{x}_{k}-\boldsymbol{x}_{j})\boldsymbol{f}_{kj}(t)m_{j}+\boldsymbol{u}(\boldsymbol{x}_{k}(t)).
\]
The linear background velocity field appears in the error representation
formula and adjoint equation. The ODE operator is augmented from $\boldsymbol{\mathcal{S}}_{\delta}^{h}[\boldsymbol{F}[\boldsymbol{x}],\boldsymbol{x}]$
to become 
\[
\boldsymbol{F}(t,\boldsymbol{x})=\boldsymbol{\mathcal{S}}_{\epsilon}^{h}[\boldsymbol{x}]+\boldsymbol{u}(\boldsymbol{x}).
\]
The addition of $\boldsymbol{u}(\boldsymbol{x})$ is then carried
through all of the error components Equations (\ref{eq:ResidualError})-(\ref{eq:QuadratureError}).
In the adjoint equation, we obtain an additional term of the form
\[
\boldsymbol{C}^{T}\cdot\boldsymbol{\phi}.
\]
For a general, nonlinear velocity field, the resulting term in the
linearized adjoint equation would depend on $\boldsymbol{x}$, e.g.
$\boldsymbol{C}^{T}(\boldsymbol{x})\cdot\boldsymbol{\phi}$. Various
choices for $\boldsymbol{C}$ lead to shear flow, straining flow,
and rotational flow \citep{pozrikidis1992boundary}.

Alternatively, it is also possible to model situations where flow
is driven by forces on the bacteria, e.g. gravitational settling of
particles. In this case, the force on each particle would be of the
form 
\[
\boldsymbol{F}[\boldsymbol{x}]+\boldsymbol{g}.
\]
Since gravitational force is independent of position and time, the
adjoint equation is not changed, and the error formulas are only modified
through a different force relation. 

Other examples of interest include cases where the forces depend on
position, or when the positions of a subset of the particles is predetermined.
The former case may occur with charged particles in an external field,
and the latter may occur if some of the particles are assumed to be
adhered to a moving boundary. Furthermore, the methods discussed here
are directly applicable kernels aside from the free-space Green's
function. For instance, the techniques developed here may be applied
to half-plane flows, or flows in a sphere where Green's functions
may be obtained through the method of images \citep{pozrikidis1992boundary}. 

In Figure \ref{fig:a,b:-Explicit-error} the accumulation of explicit
error, and the magnitude of the components of $\boldsymbol{\phi}$
are plotted versus time. We observe that the magnitude of $\boldsymbol{\phi}(t)$
grows approximately linearly over time and that the explicit error
accumulation seems to grow most rapidly at the beginning of the simulation. 

\begin{figure}
\begin{centering}
\subfloat[]{\begin{centering}
\includegraphics[width=0.3\columnwidth]{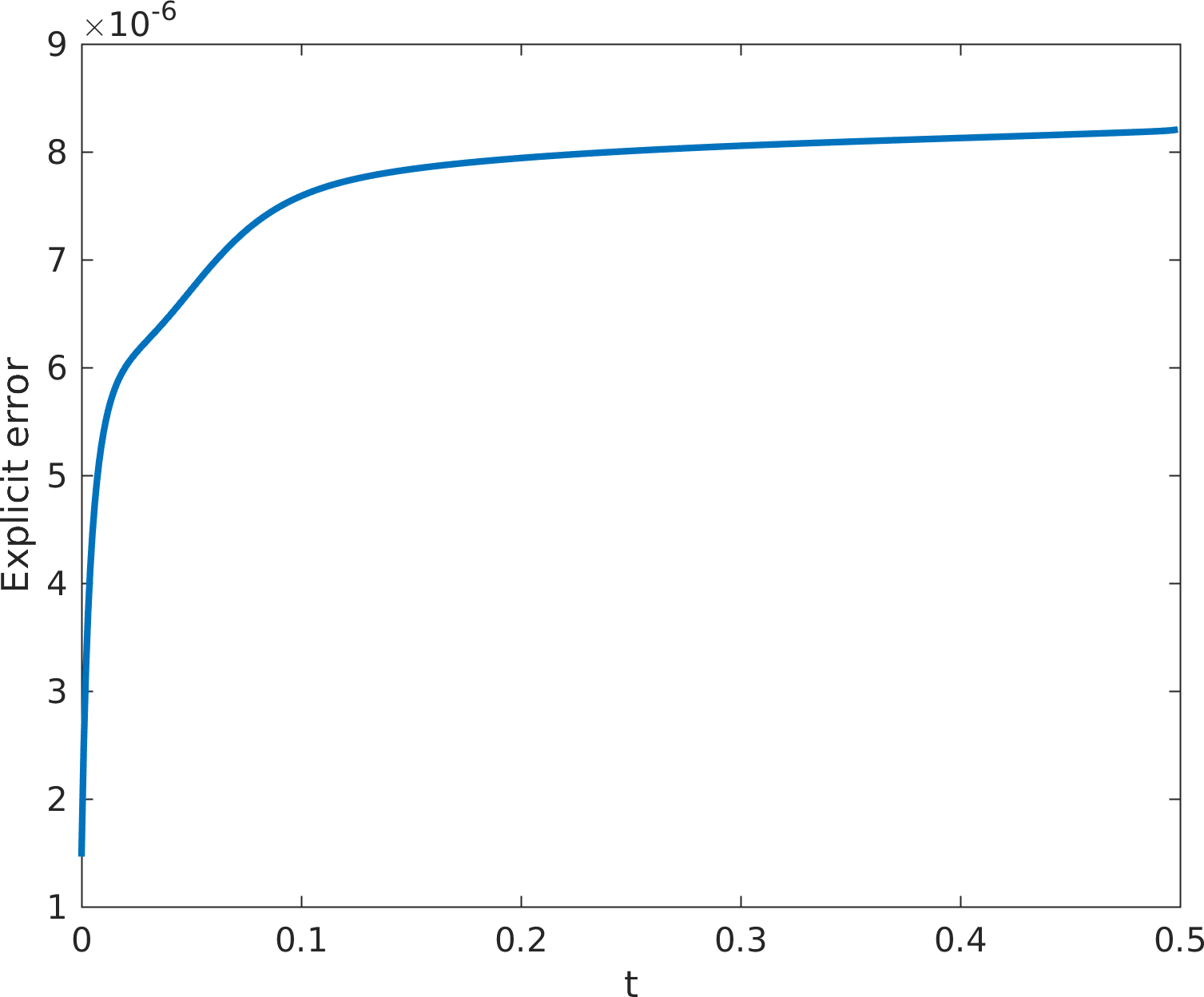}
\par\end{centering}
}\enskip{}\subfloat[]{\begin{centering}
\includegraphics[width=0.3\columnwidth]{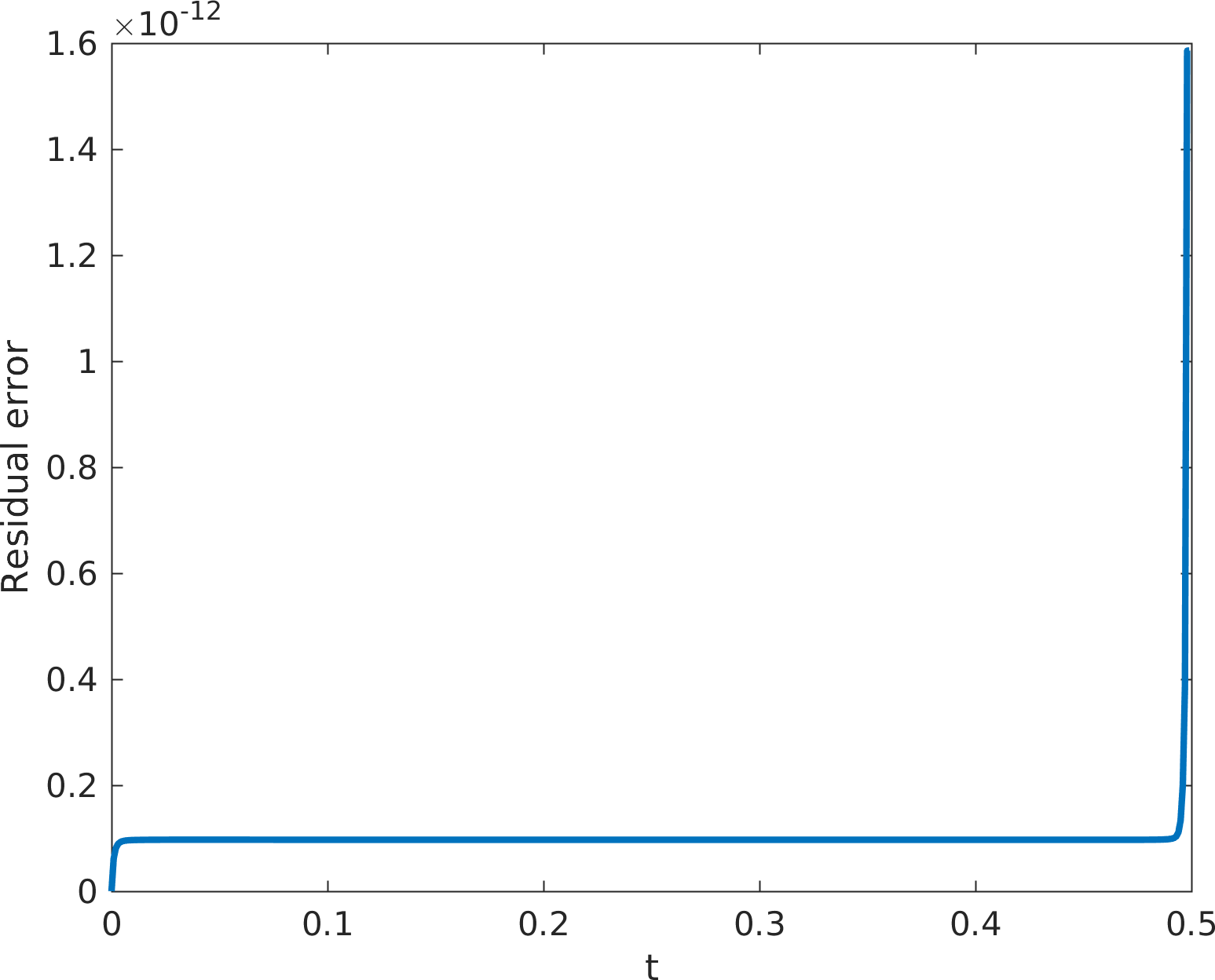}
\par\end{centering}
}\\
\subfloat[]{\includegraphics[width=0.3\columnwidth]{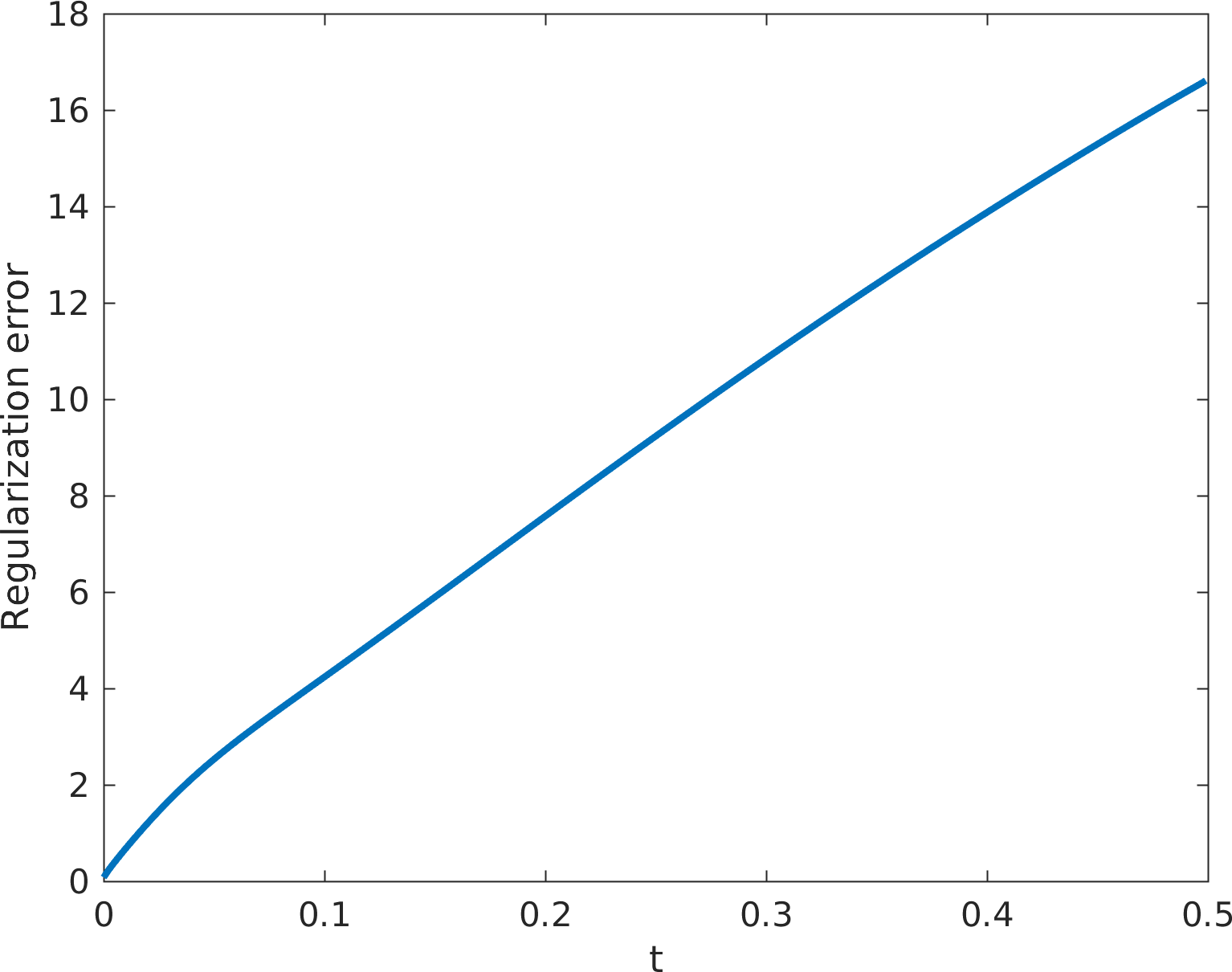}

}\enskip{}\subfloat[]{\begin{centering}
\includegraphics[width=0.3\columnwidth]{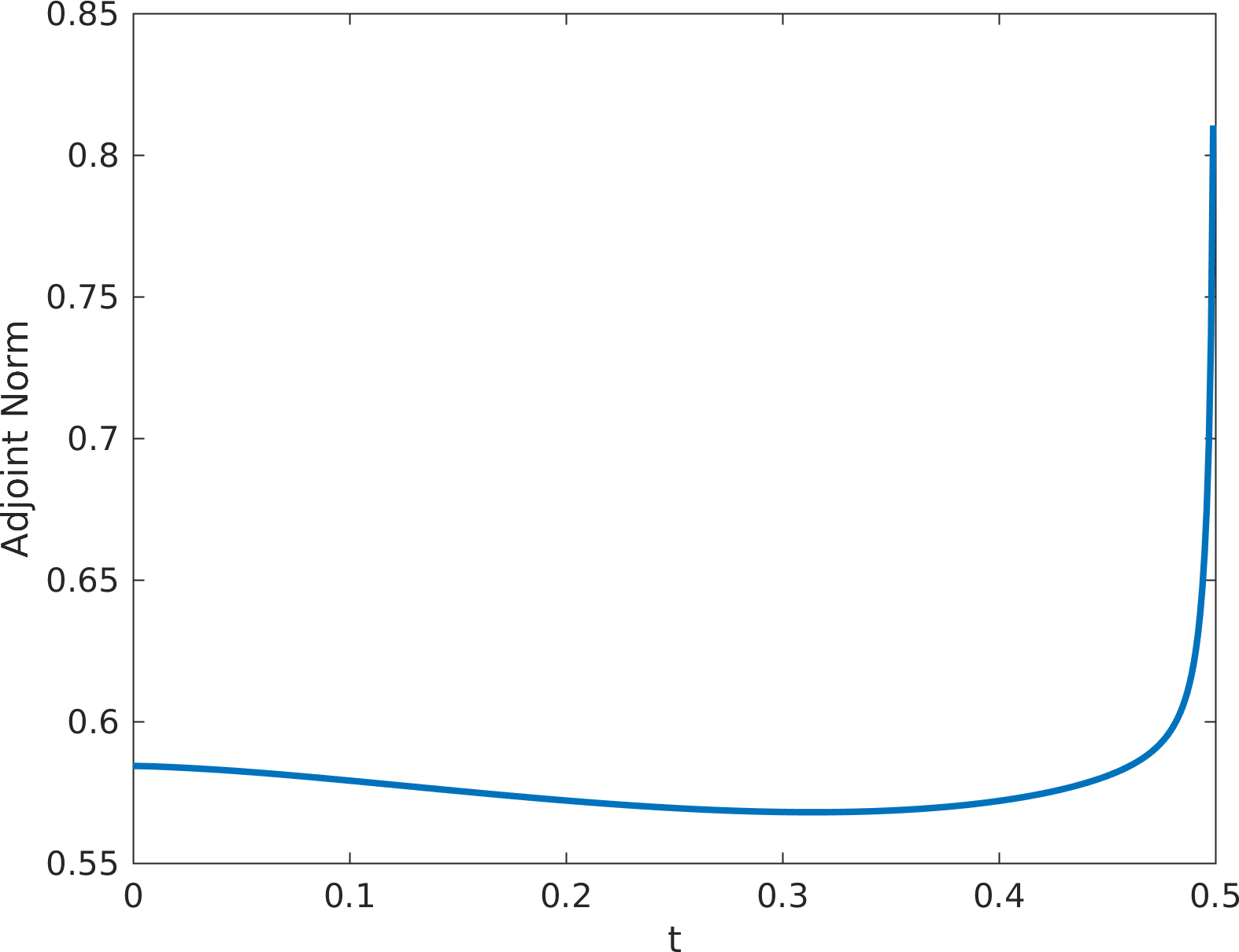}
\par\end{centering}
}
\par\end{centering}
\caption{\label{fig:a,b:-Explicit-error}Explicit Error, Residual Error, Regularization
Error, and the norm of $\boldsymbol{z}(\cdot,t)$ as functions of
time for a network of elastically connected particles.}
\end{figure}

\section{Discussion}

We have implemented the \emph{a posteriori} error estimation techniques,
originally developed in \citep{collins_posteriori_2015} in application
to a spatially discretized integrodifferential equation. Such equations
are relevant to various models in biological fluid dynamics. In particular,
the Method of Regularized Stokeslets is a popular technique for simulated
low-speed small length scale biological fluid-structure interactions.
Although the method is widely used due to its accuracy, ease of implementation,
and the fact that results from the MRS can often be compared to experimental
results, little theory exists about the well-posedness, and error
accumulation in the case of dynamically moving boundaries. Furthermore,
we are not aware of any previous studies that have looked at error
control in conjunction with the MRS. 

One ambiguity in this work is on the choice of method for extrapolating
a finite difference solution to a continuous function. We believe
there exists potential for improvements if extrapolations that are
of the same order of accuracy as the numerical solution are used instead
of lower order approximants. This seems especially likely for higher
order Runge-Kutta methods where the extrapolation leads to second
order accuracy in general, but superconvergence at the time nodes.
One way obtain such approximants may be to use a spectral discretization
in time, similar to those used in spectral deferred correction methods.
Such discretizations lead to natural choices for extrapolants that
are accurate up to the order of the discretization. We also note that
using an extrapolant that agrees at the time nodes, $t_{n}$, but
does not satisfy the property that its projection leads to $\boldsymbol{k}_{i}$
of Equation (\ref{eq:RungeKuttaMethod}) might allow for choices of
polynomial that exhibit higher order accuracy than those currently
employed by the \emph{neFEM }scheme.

Another interesting direction that may be helpful for practioners
hoping to applied the methods discussed here would be an analysis
of the impact of numerical errors in the adjoint equation solution
on the estimators. An initial investigation on this topic was conducted
in \citep{collins_posteriori_2015}, but there is certainly room for
further developments. For instance, it is agreed upon that a CFL condition
exists for the explicit discretization of the motion of an elastic
boundary in a fluid. However, it is not clear that the \emph{a posteriori}
methods we have used can necessarily capture the moments when a simulation
becomes unstable. Instability should be accompanied by rapid growth
of the adjoint solution however, it is possible that the nonlinear
behavior of the operators involved may fail to be captured by the
linearized adjoint used to obtain stability factors.

A natural next step is to combine the spatial and temporal analysis
to develop \emph{a posteriori }estimation algorithms for MRS simulations
of an elastic surface. The further extension to viscoelastic surfaces
is of interest in many applications, but poses additional, nontrivial
challenges due to the more complicated dynamic behavior of viscoelastic
materials. We are pursuing these ideas in a follow-up paper. As a
further application, we also hope to extend the methods here to quantify
the accuracy of biofilm simulations such as those described in \citep{stotsky_variable_2016}
and \citep{hammond_variable_2014}. 

\bibliographystyle{plain}
\bibliography{Bibliography}

\begin{thebibliography}{10}

\bibitem{anderson_vortex_1985}
Christopher Anderson and Claude Greengard.
\newblock On {Vortex} {Methods}.
\newblock {\em SIAM Journal on Numerical Analysis}, 22(3):413--440, June 1985.

\bibitem{aranda2015model}
Vivian Aranda, Ricardo Cortez, and Lisa Fauci.
\newblock A model of stokesian peristalsis and vesicle transport in a
  three-dimensional closed cavity.
\newblock {\em Journal of biomechanics}, 48(9):1631--1638, 2015.

\bibitem{butcher2009fifth}
JC~Butcher.
\newblock On fifth and sixth order explicit runge-kutta methods: order
  conditions and order barriers.
\newblock {\em Canadian Applied Mathematics Quarterly}, 17(3):433--445, 2009.

\bibitem{cao_posteriori_2004}
Yang Cao and Linda Petzold.
\newblock A {Posteriori} {Error} {Estimation} and {Global} {Error} {Control}
  for {Ordinary} {Differential} {Equations} by the {Adjoint} {Method}.
\newblock {\em SIAM Journal on Scientific Computing}, 26(2):359--374, January
  2004.

\bibitem{collins_posteriori_2015}
J.~B. Collins, D.~Estep, and S.~Tavener.
\newblock A posteriori error analysis for finite element methods with
  projection operators as applied to explicit time integration techniques.
\newblock {\em BIT Numerical Mathematics}, 55(4):1017--1042, December 2015.

\bibitem{cortez_accuracy_1998}
Ricardo Cortez.
\newblock On the {Accuracy} of {Impulse} {Methods} for {Fluid} {Flow}.
\newblock {\em SIAM Journal on Scientific Computing}, 19(4):1290--1302, July
  1998.

\bibitem{cortez_method_2001}
Ricardo Cortez.
\newblock The {Method} of {Regularized} {Stokeslets}.
\newblock {\em SIAM Journal on Scientific Computing}, 23(4):1204--1225, January
  2001.

\bibitem{cortez_method_2005}
Ricardo Cortez, Lisa Fauci, and Alexei Medovikov.
\newblock The method of regularized {Stokeslets} in three dimensions:
  {Analysis}, validation, and application to helical swimming.
\newblock {\em Physics of Fluids}, 17(3):031504, March 2005.

\bibitem{estep_posteriori_1995}
Donald Estep.
\newblock A {Posteriori} {Error} {Bounds} and {Global} {Error} {Control} for
  {Approximation} of {Ordinary} {Differential} {Equations}.
\newblock {\em SIAM Journal on Numerical Analysis}, 32(1):1--48, February 1995.

\bibitem{hammond_variable_2014}
Jason~F. Hammond, Elizabeth~J. Stewart, John~G. Younger, Michael~J. Solomon,
  and David~M. Bortz.
\newblock Variable {Viscosity} and {Density} {Biofilm} {Simulations} using an
  {Immersed} {Boundary} {Method}, {Part} {I}: {Numerical} {Scheme} and
  {Convergence} {Results}.
\newblock {\em CMES}, 98(3):295--340, 2014.

\bibitem{kim_microhydrodynamics:_1991}
Sangtae Kim and Seppo~J. Karrila.
\newblock {\em Microhydrodynamics: principles and selected applications}.
\newblock Butterworth-{Heinemann} series in chemical engineering.
  Butterworth-Heinemann, Boston, 1991.

\bibitem{kroese2013spatial}
Dirk~P Kroese and Zdravko~I Botev.
\newblock Spatial process generation.
\newblock {\em arXiv preprint arXiv:1308.0399}, 2013.

\bibitem{ladas1972differential}
Gerasimos~E Ladas and Vangipuram Lakshmikantham.
\newblock {\em Differential equations in abstract spaces}.
\newblock Elsevier, 1972.

\bibitem{ladyzhenskaya1969mathematical}
Olga~A Ladyzhenskaya.
\newblock {\em The mathematical theory of viscous incompressible flow},
  volume~12.
\newblock Gordon \& Breach New York, 1969.

\bibitem{lin_solvability_2017}
Fanghua Lin and Jiajun Tong.
\newblock Solvability of the {Stokes} {Immersed} {Boundary} {Problem} in {Two}
  {Dimensions}.
\newblock {\em arXiv}, 1703.03124, 2017.

\bibitem{majda_vorticity_2002}
Andrew Majda and Andrea~L. Bertozzi.
\newblock {\em Vorticity and incompressible flow}.
\newblock Cambridge texts in applied mathematics. Cambridge University Press,
  Cambridge ; New York, 2002.

\bibitem{mori_well-posedness_2017}
Y.~Mori, Analise Rodenberg, and Daniel Spirn.
\newblock Well-posedness and global behavior of the {Peskin} problem of an
  immersed elastic filament in {Stokes} flow.
\newblock {\em arXiv}, 1704.08392, 2017.

\bibitem{pozrikidis1992boundary}
Constantine Pozrikidis.
\newblock {\em Boundary integral and singularity methods for linearized viscous
  flow}.
\newblock Cambridge University Press, 1992.

\bibitem{stotsky_variable_2016}
Jay~A. Stotsky, Jason~F. Hammond, Leonid Pavlovsky, Elizabeth~J. Stewart,
  John~G. Younger, Michael~J. Solomon, and David~M. Bortz.
\newblock Variable viscosity and density biofilm simulations using an immersed
  boundary method, part {II}: {Experimental} validation and the heterogeneous
  rheology-{IBM}.
\newblock {\em Journal of Computational Physics}, 317:204--222, July 2016.

\bibitem{wrobel_modeling_2014}
Jacek~K. Wrobel, Ricardo Cortez, and Lisa Fauci.
\newblock Modeling viscoelastic networks in {Stokes} flow.
\newblock {\em Physics of Fluids}, 26(11):113102, November 2014.

\bibitem{wrobel_enhanced_2016}
Jacek~K. Wrobel, Sabrina Lynch, Aaron Barrett, Lisa Fauci, and Ricardo Cortez.
\newblock Enhanced flagellar swimming through a compliant viscoelastic network
  in {Stokes} flow.
\newblock {\em Journal of Fluid Mechanics}, 792:775--797, April 2016.

\end{thebibliography}

\end{document}